\documentclass[12pt,a4paper]{scrartcl}
\usepackage{a4wide}
\usepackage[UKenglish]{babel}
\usepackage{latexsym}
\usepackage{amssymb}
\usepackage{mathrsfs}
\usepackage{graphicx}
\usepackage[latin1]{inputenc}
\usepackage[longnamesfirst,sort&compress,comma,square,numbers]{natbib}
\usepackage{amsmath}
\usepackage{amsfonts}
\usepackage{amsthm}
\usepackage{url}
\usepackage{mathptmx}
\usepackage{helvet}
\date{}

\title{CARMA Processes driven by Non-Gaussian Noise}

\author{Robert Stelzer\thanks{Institute of Mathematical Finance, Ulm University, Helmholtzstraße 18, D-89069  Ulm, Germany. \emph{Email:} \href{mailto:robert.stelzer@uni-ulm.de}{\texttt{robert.stelzer@uni-ulm.de}}, \url{http://www.uni-ulm.de/mawi/finmath}}}

\newcommand{\bbn}{\mathbb{N}}
\newcommand{\bbz}{\mathbb{Z}}
\newcommand{\bbr}{\mathbb{R}}
\newcommand{\bbc}{\mathbb{C}}

\newcommand{\limp}{\stackrel{P}{\rightarrow}}

\newcommand{\eqd}{\stackrel{\mathscr{D}}{=}}

\newcommand{\var}{{\mathrm{var}}}
\newcommand{\cov}{{\mathrm{cov}}}

\newcommand{\bth}{\boldsymbol{\vartheta}}
\newcommand{\Y}{\mathbf{Y}}\newcommand{\y}{\mathbf{y}}
\newcommand{\X}{\mathbf{X}}
\newcommand{\Lb}{\mathbf{L}}
\newcommand{\R}{\mathbb{R}}

\newtheorem{Theorem}{Theorem}[section]

\newtheorem{Proposition}[Theorem]{Proposition}

\newtheorem{Definition}[Theorem]{Definition}

\newcounter{assumptionpara}
\newtheorem{assumptionpara}[assumptionpara]{Assumption}

\numberwithin{equation}{section}

\begin{document}
\maketitle

\begin{abstract}
We present an outline of the theory of certain L\'evy-driven, multivariate stochastic processes, where the processes are represented by  rational transfer functions (Continuous-time AutoRegressive Moving Average or CARMA models) and their applications in  non-Gaussian time series modelling. We discuss in detail their definition, their spectral representation, the equivalence to linear state space models and further properties like the  second order structure and the tail behaviour under a heavy-tailed input. Furthermore, we study the estimation of the parameters using quasi-maximum likelihood estimates for the auto-regressive and moving average parameters, as well as how to estimate the driving L\'evy process.
\end{abstract}

\vspace{0.5cm}

\noindent
\section{Introduction}
In many applications an observer (scientist, engineer, analyst) is confronted with series of data originating from one or more physical variables of interest over time. Thus, he has an observed (multivariate) time series and will often either be interested in removing (measurement) noise to extract the signal more clearly or in modelling the observed process, including its random components.

In both situations stochastic models may very well be appropriate. This is clear when one is mainly interested in removing noise, but when intending to model the observed value it is also very often appropriate to enrich a physical model by a random component to capture fluctuations and shortcomings of the physical model. The driving stochastic process (the ``noise'') may have interest on its own (as is the case with economic models), or it has to be modelled well to extract the interesting information as well as possible (e.g., as is common practice in telecommunication links)

The easiest way to obtain a model with randomness for the variables of interest would be to assume that all observed values are  independent and identically distributed (iid) random variables or that they follow a physical model plus iid noise. However, in most series observed consecutive values are heavily dependent and thus more sophisticated models are needed. A flexible but at the same time very tractable class of models is given by linear random processes. In the discrete time setting these models are well-known as autoregressive moving average (ARMA) processes and they are given in terms of a general order linear difference equation where an iid noisy input sequence introduces all randomness. The latter is also referred to as linear filtering of a white noise.

In many situations it is more appropriate to specify a model in continuous time rather than in discrete time. These include high-frequency data, irregularly spaced data, missing observations or situations when estimation and inference at various frequencies is to be carried out. Moreover, many physical models are formulated in continuous time and, hence, such an approach is often more natural. 

In the following we consider linear random processes in continuous time, referred to as con\-tinu\-ous time  autoregressive moving average (CARMA) processes. Intuitively, they are given as the solution to a higher order system of linear differential equations  with a stochastic process as the input, which can be seen as linearly filtering the random input.

One important  question is which random input to take in the continuous time set-up. Clearly, the random process should correspond in some sense to the idea of white noise. Understanding the latter in the strict sense means using independent increments, in the weak sense it means uncorrelated increments and so the variance has to be finite. Recall that for random variables uncorrelatedness is equivalent to independence only if the random variables are Gaussian, i.e. they have a normal distribution. A linear random process driven by Gaussian white noise has again Gaussian distributions. However, in many situations it is not appropriate to assume Gaussianity of the variables of interest, since the observed time series often exhibit features like skewness or heavy-tails (i.e. very high or low values are far more likely to occur than in the Gaussian setting), which contradict the Gaussian assumption. Demanding uncorrelated but not necessarily independent  increments does not lead to a nice class of processes nor to nice theoretical results.

Hence, a good modelling strategy where the resulting process is reasonably tractable and the driving process' probability distribution is allowed to have ``fat tails'' is to demand that the random input shall have independent as well as stationary increments, i.e. increments over time intervals of the same length have the same distribution. They then have a time homogeneity feature and resemble the iid noise of the discrete time set-up. The resulting class of possible driving processes are the so-called L\'evy processes,
which have been studied in detail and form a both highly versatile and highly tractable family. An interesting feature is that linear processes driven by general L\'evy processes may exhibit jumps and thus allow the modelling of abrupt changes, whereas Gaussian linear processes have continuous sample paths.

In the remainder of this paper we proceed as follows. First, we introduce L\'evy processes in detail. Thereafter, we give a proper definition of CARMA processes, discuss their relation to linear filtering via a stochastic Fourier (spectral) representation and summarize central properties of CARMA processes. Next, we briefly explain the equivalence to linear state space models and the relation to stochastic control and signal processing.  Finally, we discuss the statistical estimation of the parameters and the underlying L\'evy process and conclude with some additional remarks.

Throughout we will focus on developing the main ideas for CARMA processes. For more mathematical details as well as comprehensive references we refer the interested reader to the original literature especially the works \cite{Arato1982}, \cite{Brock1,Brockwell2001,Brockwell2009}, \cite{BrockwellDavisYang2007,BrockwellDavisYang2007b,BrockwellDavisYang2009}, \cite{Brockwell:Lindner:2009b}, \cite{Brockwelletal2005}, \cite{Marquardtetal2005} and \cite{Schlemm:Stelzer:2010a}. For a historic perspective the monograph \cite{Priestley1} may be interesting as well as \cite{Doob1944} which is the first paper where Gaussian CARMA processes appeared under the name of Gaussian processes with rational spectral density.

\section{L\'evy processes}

A L\'evy process $L=(L_t)_{t\in\bbr^+}$ is a stochastic process with independent and stationary increments. In the following we consider only L\'evy processes taking values in the $m$-dimensional vector space $\bbr^m$ (with $\bbr$ the real numbers and $m$ some positive integer). Note that a stochastic process $(X_t)_{t\in\bbr^+}$ can be either seen as a family of random variables indexed by the positive real numbers $\bbr^+$ or as a random function mapping the positive real numbers to $\bbr^m$. More precisely we have the following definition:
\begin{Definition}
An $\bbr^m$-valued stochastic process $L=(L_t)_{t\in\bbr^+}$ is called \emph{Lévy process}  if
\begin{itemize}
 \item
  $L_0=0$ a.s.,
 \item
  $L_{t_2}-L_{t_1}$,$L_{t_3}-L_{t_2}$,\ldots, $L_{t_n}-L_{t_{n-1}}$ are independent  for all $n\in\bbn$ and $t_1,t_2,\ldots,t_n\in\bbr+$ with $0\leq t_1< t_1<\ldots<t_n$,
 \item
  $L_{t+h}-L_t\eqd L_{s+h}-L_s$ for all $s,t,h\in\bbr^+$ (``$\eqd$'' denoting equality in distribution),
 \item
  $L$ is continuous in probability, i.e. for all $s\in\bbr^+$ we have $L_t-L_s\limp 0$ as $t\to s$.
\end{itemize}
\label{def:levy}
\end{Definition}

It can be shown
 (cf. \cite{Sato1999} for a detailed proof)
 that the class of L\'evy processes can be characterized fully at the level of ``characteristic functions'', which we now introduce. Let $<\cdot,\cdot>$ indicate the natural inner product in $\bbr^{m}$ and $X$ is an $\bbr^{m}$-valued random variable, then its {\it characteristic function} is defined as $\psi_{X}(u)=E \left(e^{i<u,X>}\right)$).  The characteristic function of a L\'evy process can always be represented  in the {\it L\'evy-Khintchine form}
\begin{align}
 E\left(e^{i\langle u,{L}_t\rangle}\right)&=\exp\{t\psi_L(u)\},\;\;\forall \,t\geq 0,\, u\in\bbr^m,\\
\intertext{ with }
\psi_L(u)&=i\langle\gamma,u\rangle-\frac{1}{2}\langle u,\Sigma_G u\rangle+\int\limits_{\bbr^m}
(e^{i\langle u,x\rangle}-1-i\langle u,x\rangle 1_{[0,1]}(\|x\|)\,\nu(dx),
\end{align}
where $\gamma\in\bbr^m$, $\Sigma_G$ is a $m\times m$ positive semi-definite matrix and $\nu$ is a measure on $\bbr^m$ that satisfies
$\nu(\{0\})=0$ and  $\int\limits_{\bbr^m}(\|x\|^2\wedge 1)\,
\nu(dx)<\infty.$
The measure $\nu$ is referred to as the L\'evy measure of ${L}$ and $\|x\|^2\wedge 1$ is short for $\min\{\|x\|^2,1\}$. Finally, $1_A(x)$ generically denotes the indicator function of a set $A$, i.e. the function which is one if $x$ is an element of $A$ and zero otherwise. Together $(\gamma,\Sigma_G,\nu)$ are referred to as the characteristic triplet of $L$.

Regarding the paths of a L\'evy process, i.e. the ``curve of $L$ as a function of time $t$, it can be shown that without loss of generality, a L\'evy process may be assumed to be right continuous and have left limits.

It should be noted that many well-known stochastic processes are L\'evy processes. Examples are Brownian motion, also referred to as the Wiener process or ``Gaussian white noise'', the Poisson process, which has jumps of size one and remains constant in between the jumps, which occur after iid exponentially distributed waiting times, and $\alpha$-stable L\'evy motions, sometimes called L\'evy flights. Compound Poisson processes are Poisson processes where the fixed jump size one is replaced by random iid jump sizes independent of the interarrival times of the jumps. It can be shown that all L\'evy processes arise as limits of such compound Poisson processes. 

\begin{figure}[tp]
  \centering
  \includegraphics[width=0.48\textwidth]{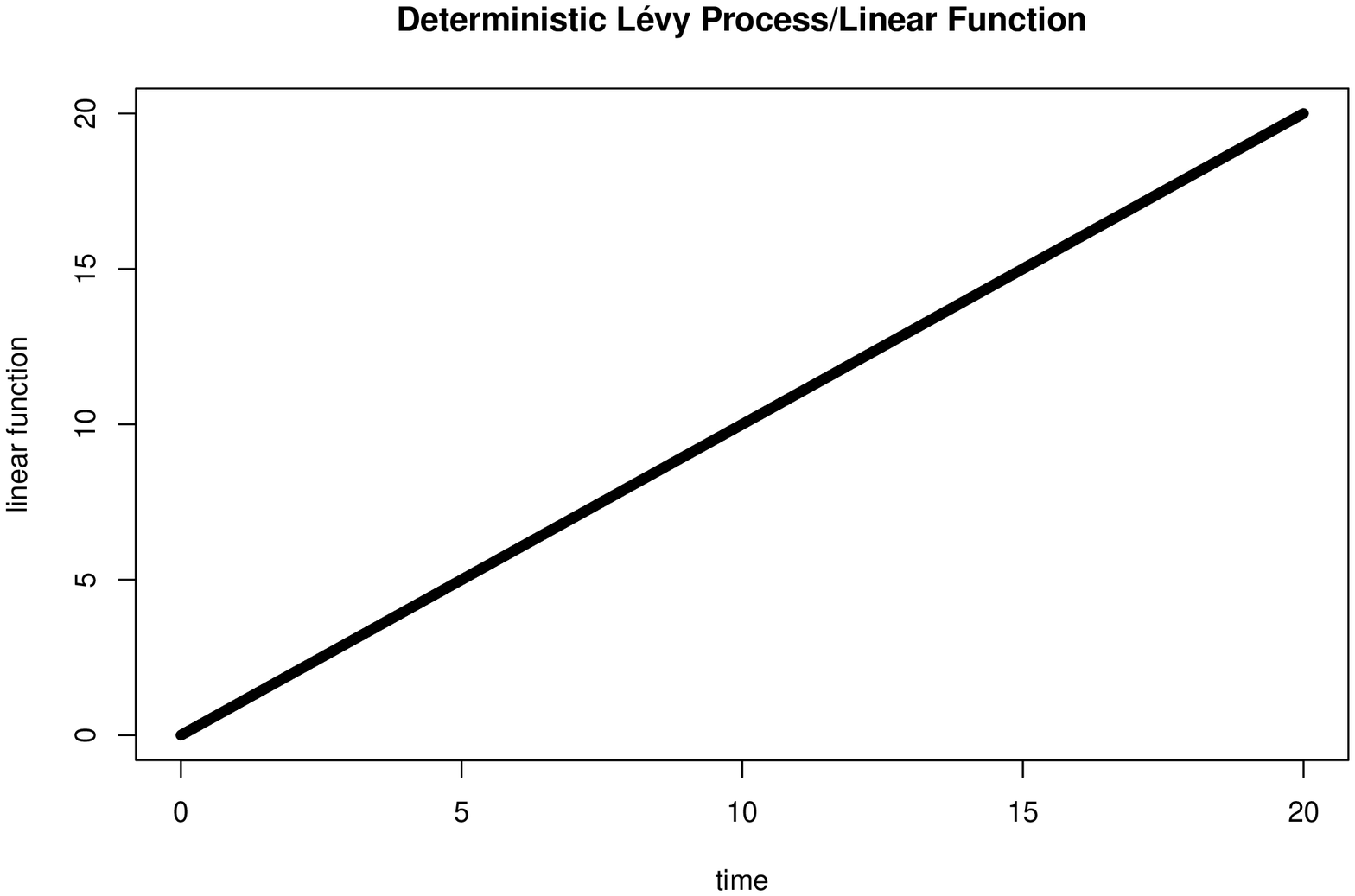}
\includegraphics[width=0.48\textwidth]{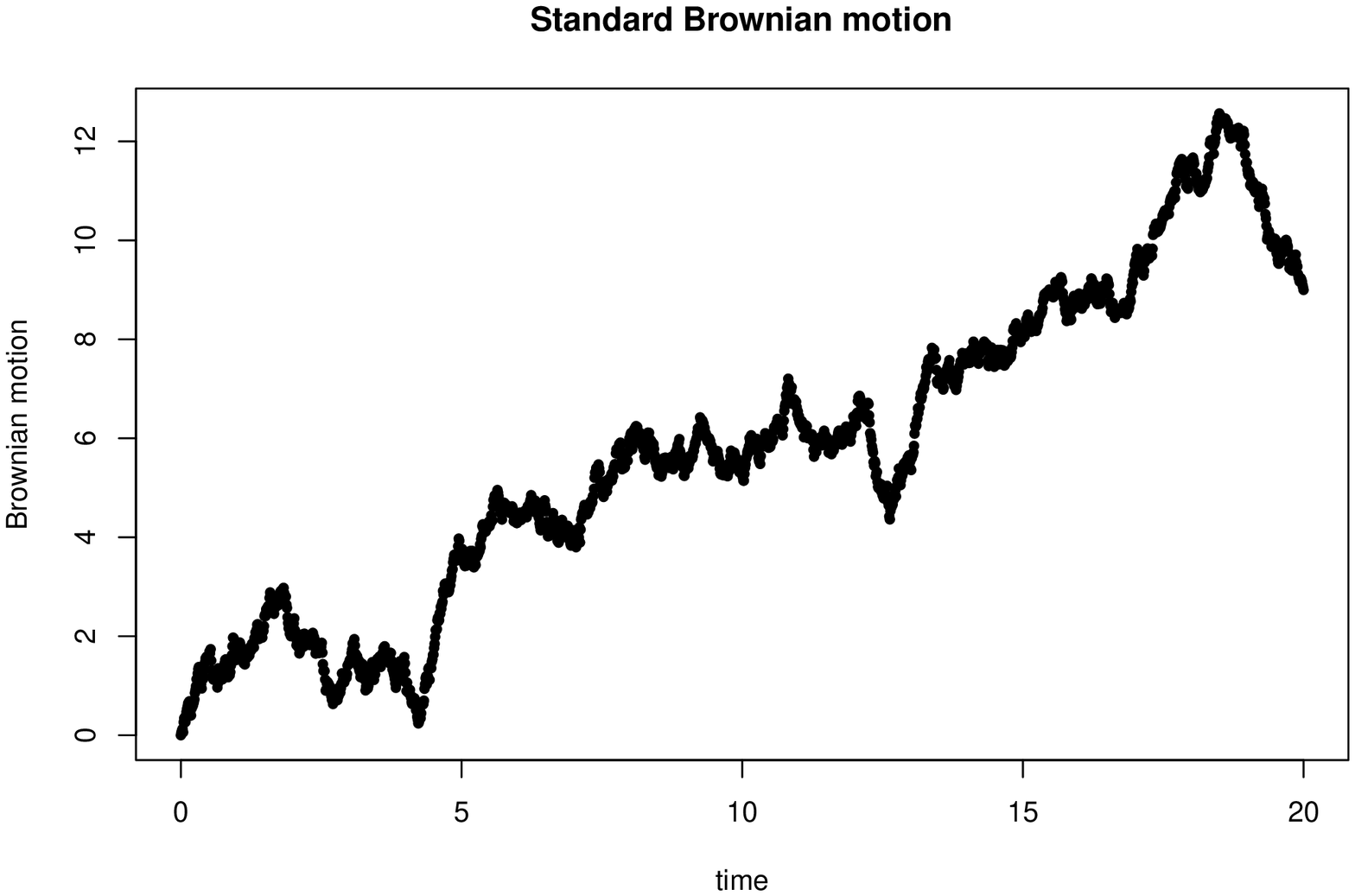}
\includegraphics[width=0.48\textwidth]{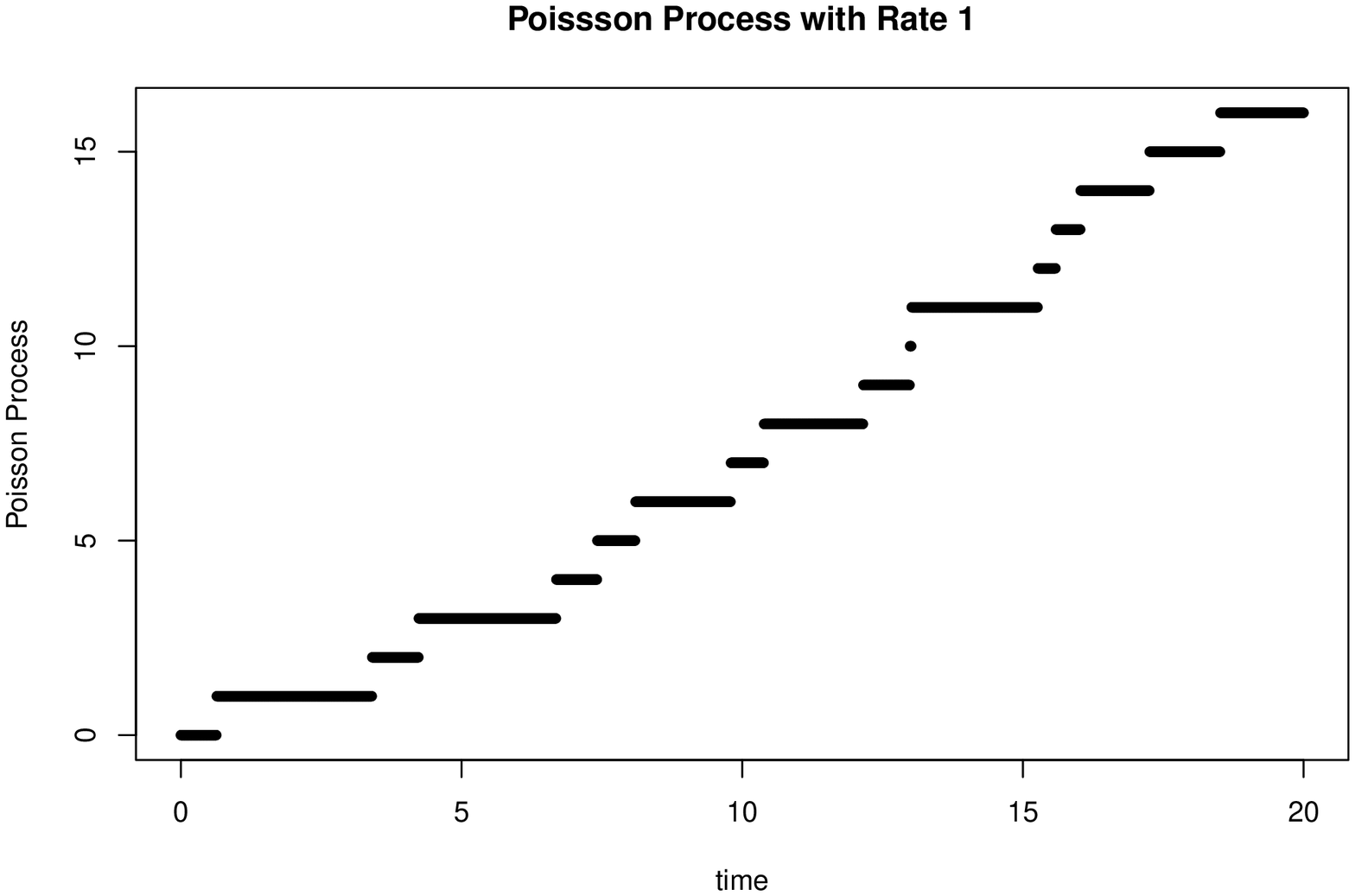}
 \includegraphics[width=0.48\textwidth]{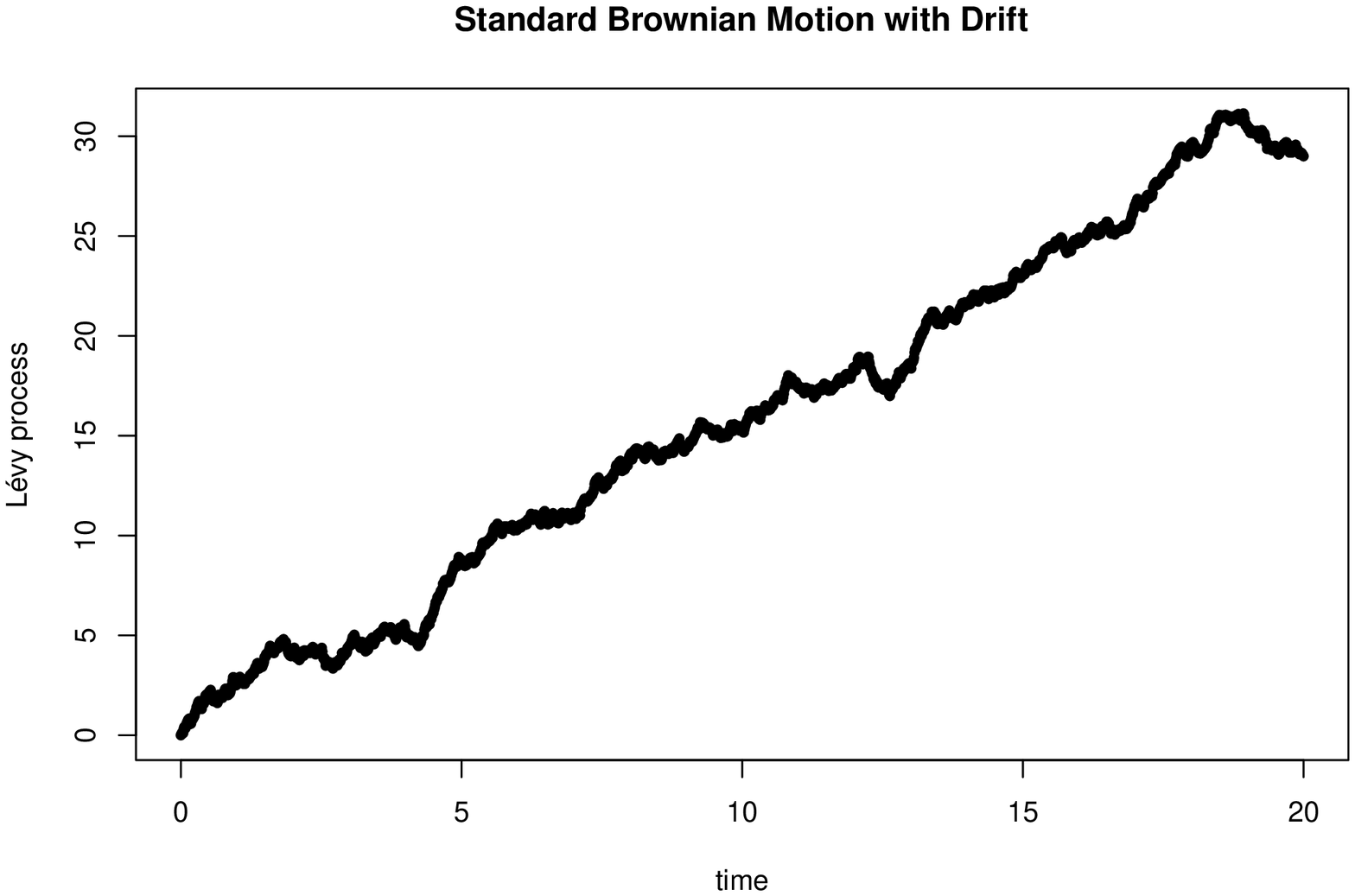}
 \includegraphics[width=0.48\textwidth]{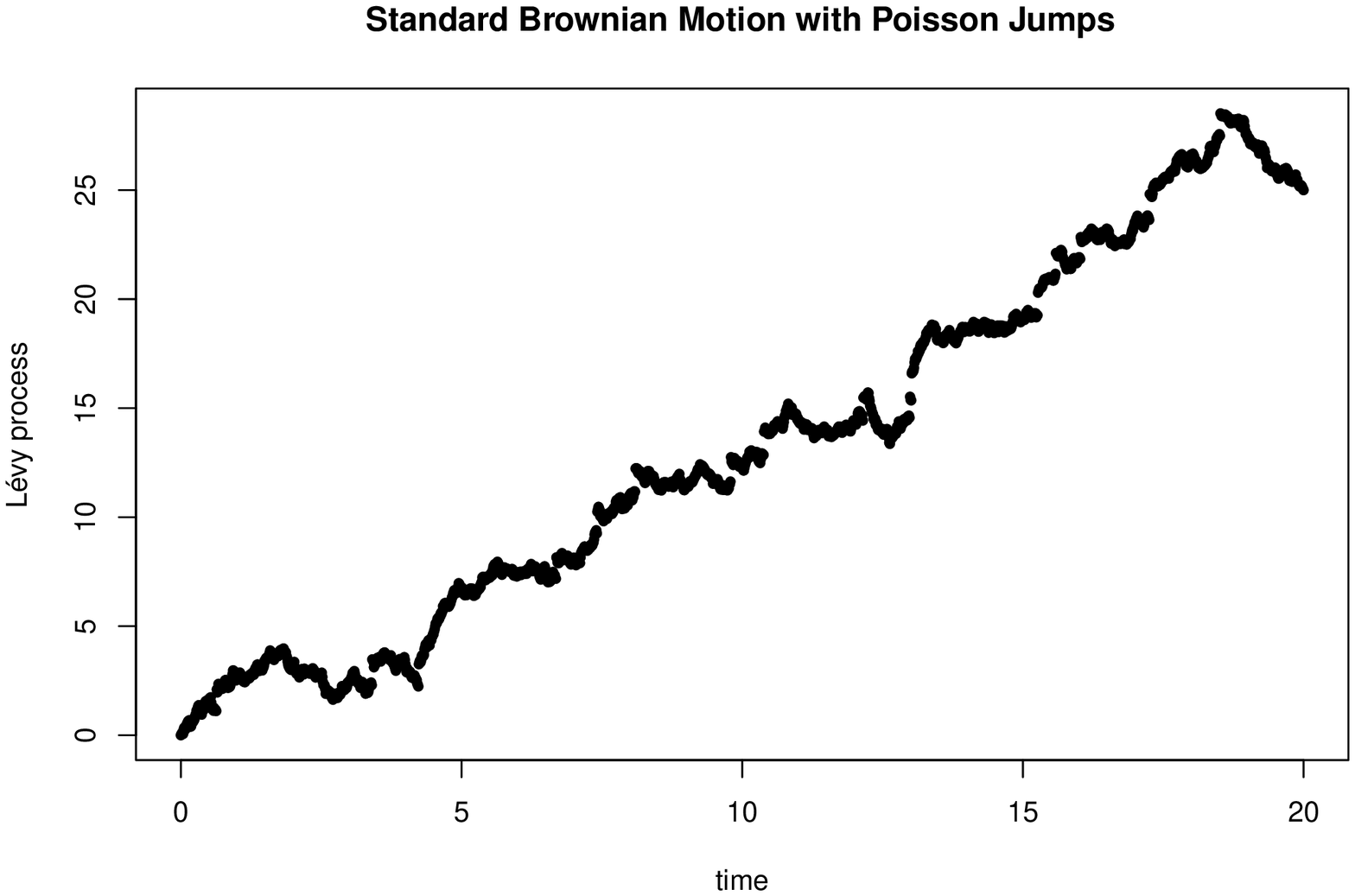}
\includegraphics[width=0.48\textwidth]{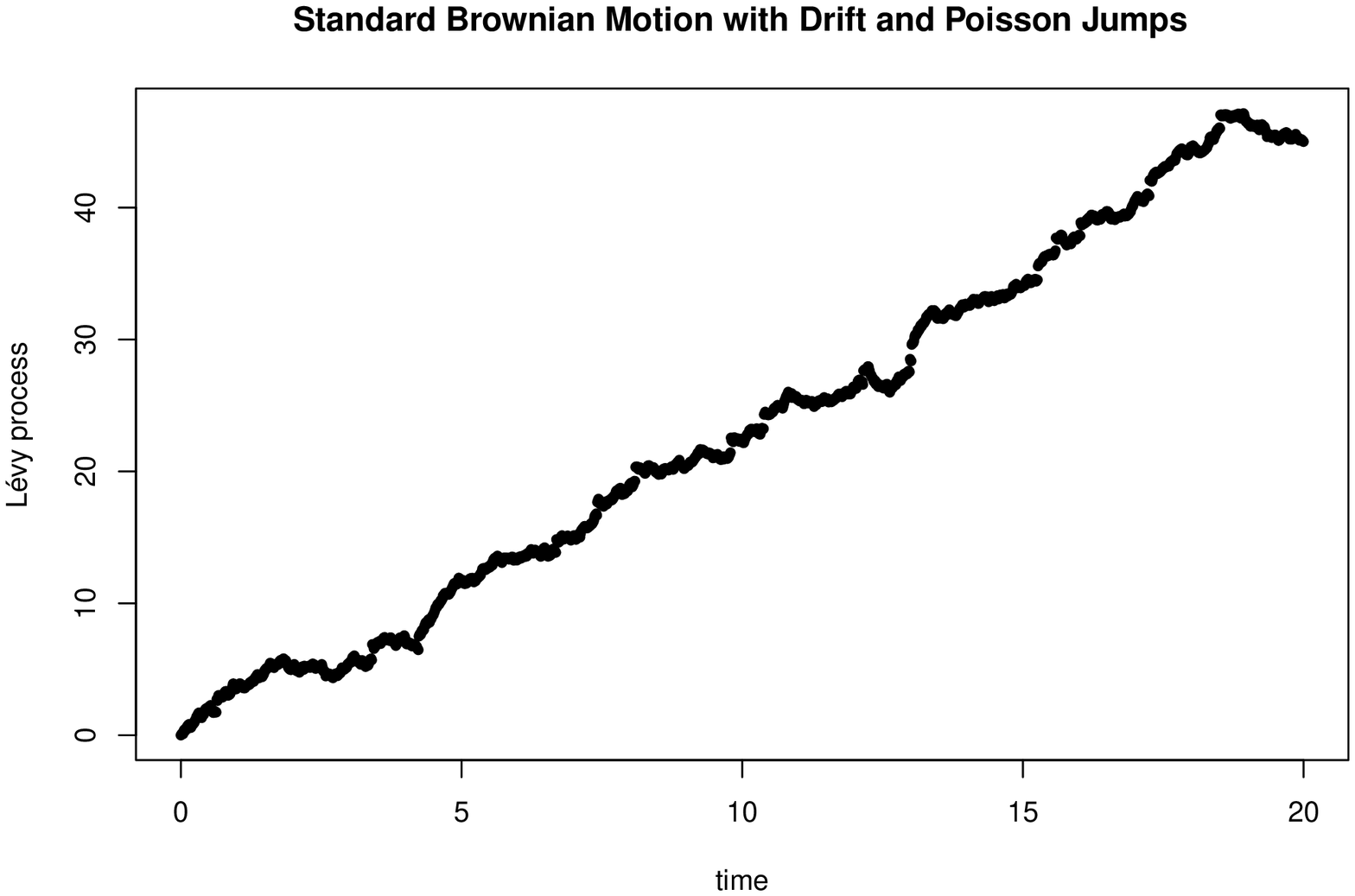}
  \caption{A L\'evy process and its components: The complete L\'evy process is depicted in the lower right display. In the upper row the deterministic drift component is depicted on the left and the standard Brownian motion component on the right. The left display in the middle row shows the standard (rate one) Poisson component and the right one the Brownian motion and the deterministic component added together. In the last row on the left the Brownian component plus the Poisson jumps are depicted.\newline
  Note that the scaling of the $y$-axis is different in the individual plots.}
  \label{fig:Levycomp}
 \end{figure}

A better understanding of what L\'evy processes really are is provided by the L\'evy-It\^o decomposition of their paths. It states that a L\'evy process is the sum of the deterministic linear function $\gamma t$, a Brownian motion with covariance matrix $\Sigma_G$, the sum of the big jumps which form a compound Poisson process and the compensated sum
 of the small jumps (i.e. the sum of the small jumps minus their expected value). The quantity $\nu(A)$ gives for any measurable set $A\subset\bbr^m$ the expected number of jumps with size in $A$ occurring in a time interval of length one. In Figure  \ref{fig:Levycomp} a univariate L\'evy process which is the sum of  the linear function $t$, in this case with $\gamma=2$, a standard Brownian motion, with $\Sigma_G=1$, and a Poisson process, with $\nu(\{1\})=1, \,\nu(\bbr\backslash\{1\})=0$ is depicted together with its individual components.

Whenever $\int_{\bbr^m}(\|x\|\wedge 1)\nu(dx)<\infty$, we can replace the compensated sum of small jumps simply by the sum of the small jumps adjusting also the slope of the deterministic component. We have actually already done this in Figure  \ref{fig:Levycomp} where the resulting slope of the deterministic function is $\gamma-\int_\bbr x\nu(dx)=1$. If $\nu(\bbr)<\infty$, we have finitely many jumps in any bounded time interval and the jumps form actually a compound Poisson process. Otherwise, we have infinitely, but countably many jumps in any bounded time interval. The reason why we have in general a component referred to as ``the compensated sum of the jumps'' (i.e., it results from a certain limiting procedure see e.g. \cite{Sato1999}) is that in general the jumps are not summable. This is equivalent to the fact that the paths have infinite variation, like Brownian motion. Infinite variation intuitively means that the curve described by the stochastic process over finite time intervals has an infinite length. Clearly, this means that the fluctuations of the process over small time intervals are rather vivid.

\begin{figure}[tp]
 \centering
 \includegraphics[width=\textwidth,bb=0 0 1049 692]{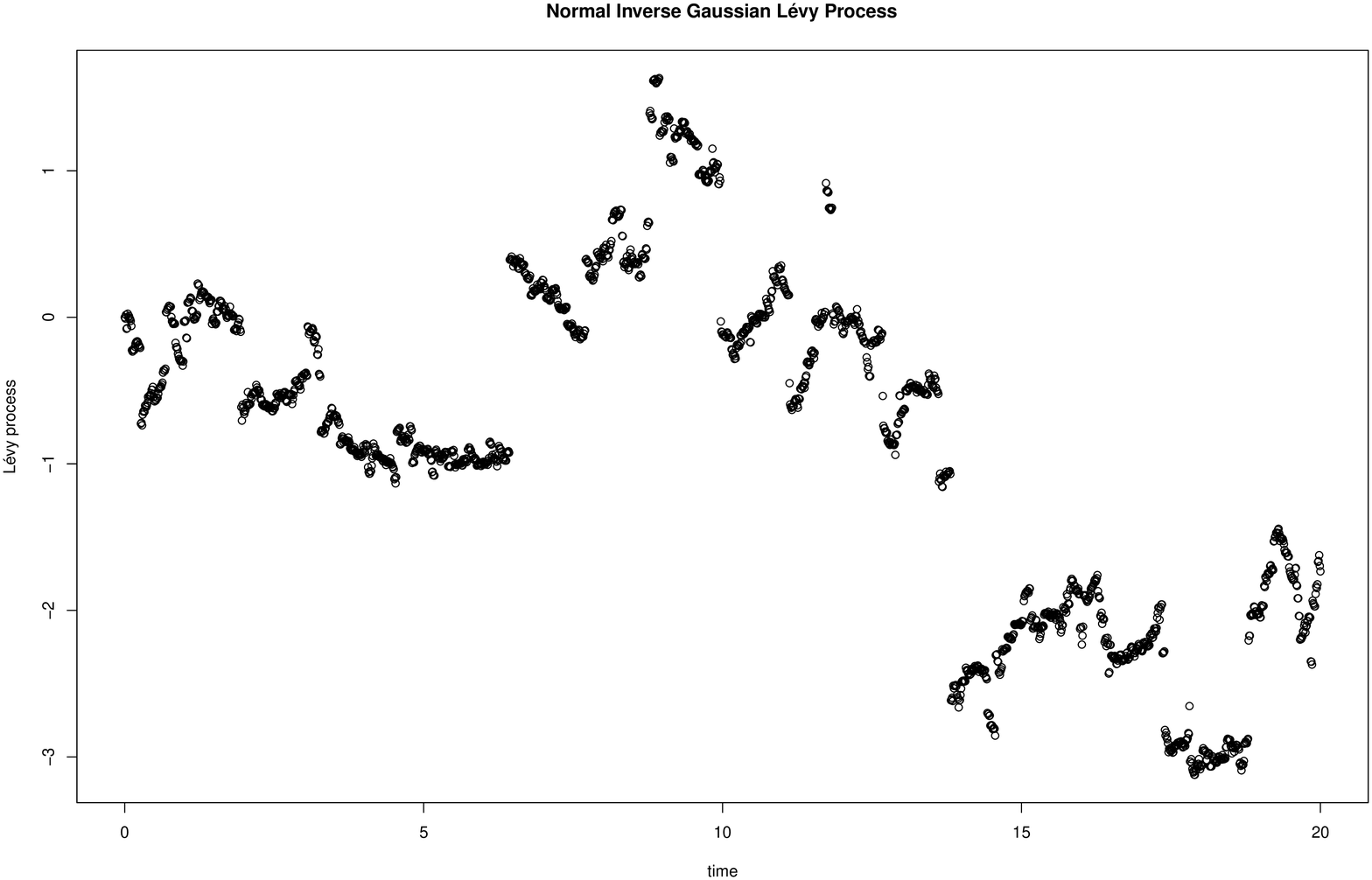}
\caption{Simulation of a Normal Inverse Gaussian (NIG) L\'evy process.}
 \label{fig:nig}
\end{figure}

\begin{figure}[tp]
\includegraphics[width=\textwidth,bb=0 0 1049 692]{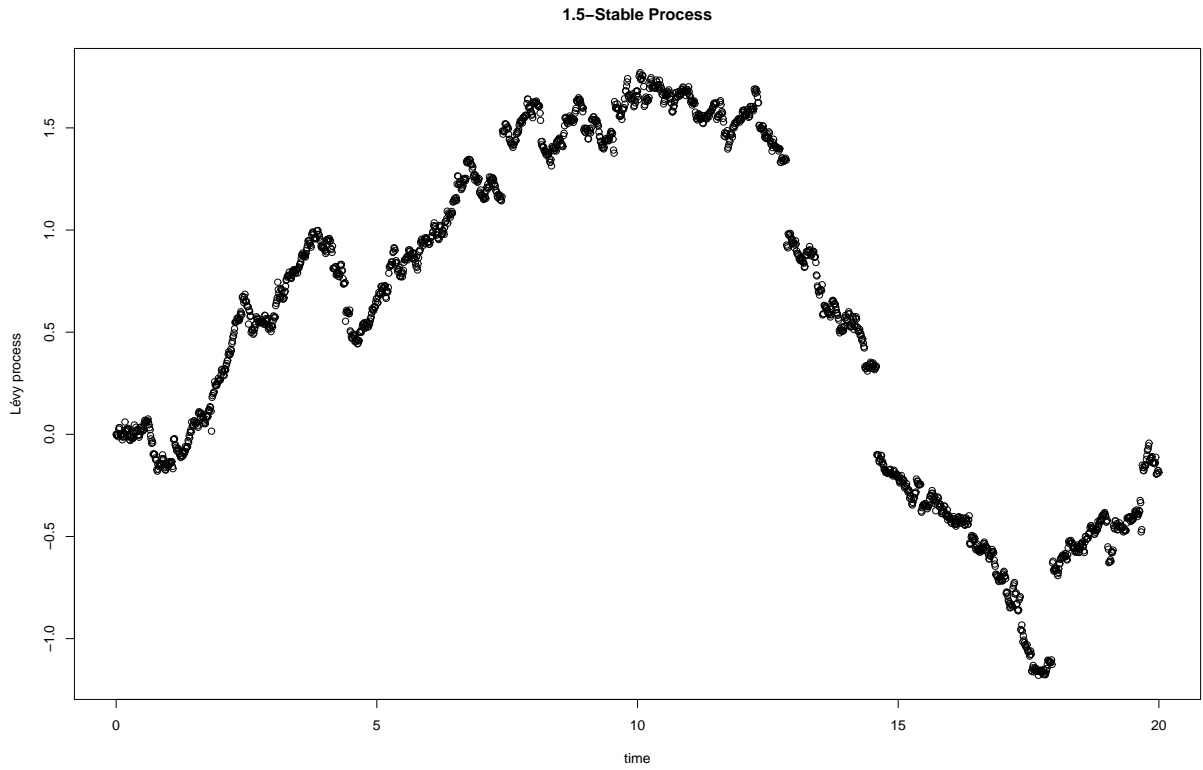}

 \includegraphics[width=\textwidth,bb=0 0 1049 692]{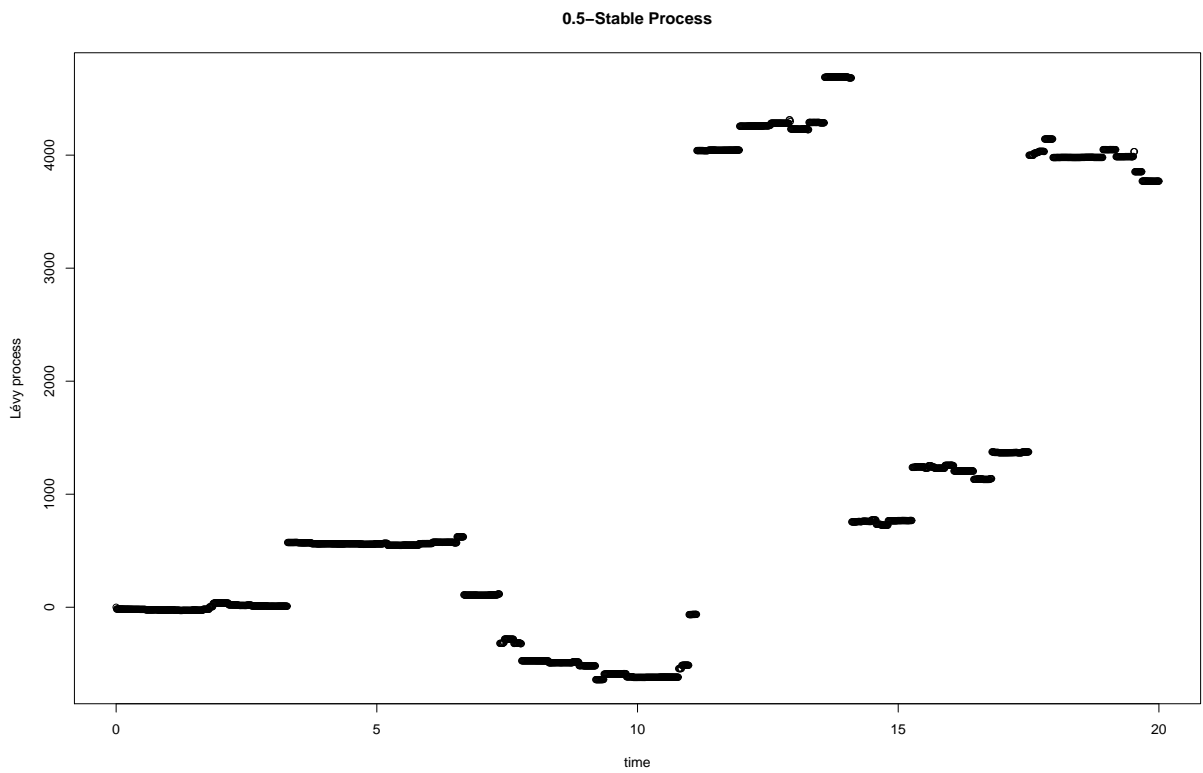}
\caption{Simulations of stable L\'evy processes.  A 1.5-stable L\'evy process is depicted in the upper row and a 0.5-stable in the lower one.}
 \label{fig:stable}
\end{figure}

In Figures \ref{fig:nig} and \ref{fig:stable} you can see simulations of different pure jump L\'evy processes, i.e. in these cases $\gamma=0$ and $\Sigma=0$. So there is neither a deterministic drift nor a Brownian motion present. All these processes have infinite activity, i.e. infinitely many jumps in any time interval. Figure \ref{fig:nig} depicts a so-called normal inverse Gaussian L\'evy process which has heavier tails than a Brownian motion, but still is rather tame, because it has finite moments of all orders, i.e. $E(|L_t|^r)<\infty$ for all $t,r\in\bbr^+$ and also some exponential moments. In contrast to this the stable processes of Figure \ref{fig:stable} have very heavy tails, because they do not have a finite variance and the 0.5-stable processes does not even have a finite mean. Whereas the NIG and 1.5-stable processes have infinite variation, the small jumps of the 0.5-stable L\'evy process are summable.


Most of the time we will work with L\'evy processes defined on the whole real line, i.e. indexed by $\bbr$ not $\bbr^+$. They are obtained by taking two independent copies of a L\'evy process and reflecting one copy at the origin.

For detailed expositions on L\'evy processes we refer to \cite{Applebaum2004}, \cite{Bertoin1998},  \cite{Kyprianou2006} or \cite{Sato1999}.

\section{Definition of CARMA processes and spectral representation}
On the intuitive level one wants to be able to interpret a $d$-dimensional {\it CARMA($p,q$) process $Y$}  as the stationary solution to the $p$-th order linear differential equation
\begin{align}\label{eq:carmadiffeq}
P(D)Y_t&=(D^p+A_1D^{p-1}+\ldots+A_p)Y_t\\&=(B_0D^q+B_1D^{q-1}+\ldots+B_q)DL_t={Q(D)DL_t},
\end{align}
where
 the driving input $L$ is an $m$-dimensional L\'evy process,
 $D$ denotes differentiation with respect to $t$,
  and the coefficients $A_1,\ldots, A_p$ are $d\times d$ matrices and $B_0,\ldots,B_q$ are $d\times m$ matrices. The polynomials $P(z)=z^p+A_1z^{p-1}+\ldots+A_p$ and $Q(z)=B_0z^q+B_1z^{q-1}+\ldots+B_q$ with $z\in\bbc$ are referred to as the auto-regressive and moving average polynomial, respectively. Finally, $p,q\in\bbn$ are the auto-regressive and moving average order.

However, the paths of non-deterministic L\'evy processes are not differentiable and so the above equation cannot directly provide a rigorous mathematical definition. 
Let us briefly consider the case $(p,q)=(1,0)$ in which case  the resulting process is actually called an Ornstein-Uhlenbeck (OU) process. In the univariate case it is given by the differential equation
\[
D Y_t=a Y_t+ DL_t
\]
where $a$ is a real number. So what we basically want is that the change of $Y$ over an infinitesimal time interval is $a$ times the current value of the process times the ``length of the infinitesimal time interval'' plus the change of the L\'evy process over the infinitesimal time interval. Rephrasing this idea in the precise language of stochastic differential equations (see e.g. \cite{Protter2004}) we obtain
\[
dY_t=aY_t dt+dL_t.
\]
Using the theory of stochastic differential equations (SDEs) it is easy to see that this SDE has a unique solution given by 
\[
 Y_t=e^{at}Y_0+e^{at}\int_0^t e^{-as}dL_s.
\]

For general orders $(p,q)$ one could to some extent use a similar reasoning to arrive at a precise definition of CARMA processes. However, we shall take a more elegant route. First note that the differential operators on the auto-regressive side of \eqref{eq:carmadiffeq} act like integration operators on the moving average side. Hence, they offset the differential operators of the moving average side acting on the L\'evy process. Since L\'evy processes are not differentiable, we effectively have to integrate at least as often as we differentiate to be able to make sense of \eqref{eq:carmadiffeq}. Hence, a necessary condition ensuring the proper existence of CARMA processes is $p>q$.

In order to obtain a rigorous definition of CARMA processes our strategy here shall be to switch from the time domain to the frequency domain where the main tool is the following spectral representation of a L\'evy process. Here and in the following we denote by $A^*$ for a matrix (or vector) $A$ the Hermitian, i.e. the complex conjugate transposed matrix.
\begin{Theorem}[\cite{Marquardtetal2005}]\label{th:rom}
Let $(L_t)_{t\in\bbr}$ be a  \emph{square integrable $m$-dimensional L\'evy process} with mean $E[L_1]=0$
(which implies $E[L_{t}]=0$ for all $t$)
 and variance $E[L_1L_1^*]=\Sigma_L$ at $t=1$.
 Then there exists a {unique $m$-dimensional random orthogonal measure  $\Phi_L$} with spectral measure $F_L$ 
such that $E[\Phi_L(\Delta)]=0$ for any bounded Borel set $\Delta$,
$ F_L(dt)=\frac{\Sigma_L\,}{2\pi}dt\label{specmeasure}$
and 
\[{L_t=\int_{-\infty}^\infty \frac{e^{i\mu t}-1}{i\mu}\,\Phi_L(d\mu)},\: t\in\bbr.
\]
The random measure $\Phi_L$ is uniquely determined by 
\begin{equation}
 \label{def:phil}\Phi_L([a,b))=\int\limits_{-\infty}^\infty \frac{e^{-i\mu a}-e^{-i\mu b}}{2\pi i\mu}\,dL_{\mu}
\end{equation}
 for all $-\infty<a<b<\infty$.
\end{Theorem}
The random orthogonal measure $\Phi_L$ can intuitively be thought of as the ``Fourier transform'' of the L\'evy process. If $L_{t}$ is a Brownian motion, then $\Phi_L([0,t))$ is again a Brownian motion. For general L\'evy processes rather little can be said about the properties of $\Phi_L$. For example, it is known that $\Phi_L$ has second-order stationary and uncorrelated increments, but the increments are neither independent nor stationary in a strict sense, see \cite{FuchsStelzer2010}.  

In the spectral domain we can now interpret differentiation (and integration) as linear filtering noting that a formal interchange of differentiation and integration gives ``$DL_t=\int_{-\infty}^\infty {e^{i\mu t}}\,\Phi_L(d\mu)$''. It can be shown that the resulting process is well-defined whenever the linear filter is square integrable. Thus we obtain as definition for  ``$Y(t)=P(D)^{-1}Q(D)DL(t)$'':
\begin{Definition}[CARMA Process, \cite{Marquardtetal2005}]\label{def:mcarma} Let $L=(L_t)_{t\in\bbr}$ be a two-sided 
square integrable $m$-dimensional L\'evy-process with $E[L_1]=0$ and $E[L_1L_1^*]=\Sigma_L$.
A $d$-dimensional L\'evy-driven continuous time autoregressive moving average
process $(Y_t)_{t\in\bbr}$ of order ($p,q$) with  $p,q\in\bbn_0$ and  $p>q$ (\emph{CARMA($p,q$) process}) is defined as
\begin{eqnarray}
Y_t&=&\int\limits_{-\infty}^{\infty} e^{i\mu t} P(i\mu)^{-1}Q(i\mu)\,\Phi_L(d\mu),
\quad t\in\bbr, \qquad\text{where}
\label{specMCARMA}
\\
P(z):&=&I_mz^p+A_1z^{p-1}+...+A_p,\nonumber
\\
Q(z):&=&B_0z^q+B_1z^{q-1}+....+B_q\quad\text{and} \nonumber
\end{eqnarray}
$\Phi_L$ is the L\'evy orthogonal random measure of Theorem \ref{th:rom}.
Here
$A_j\in M_m(\bbr)$, $j=1,...,p$ and 
$B_j\in M_{d,m}(\bbr)$ are matrices satisfying  $B_q\not =0$ and $\mathcal{N}(P):=\{z\in\bbc:\det(P(z))=0\}\subset \bbr\backslash\{0\}+i\bbr$ (i.e. the autoregressive polynomial has no zeros on the complex axis).
\end{Definition}
Referring to the explicit construction of the random orthogonal  measure $\Phi_L$, one can easily show that the above defined CARMA processes are necessarily stationary (in the strict sense, i.e. the distributions are left unchanged by a time shift). Since by construction any CARMA process in the sense of Definition \ref{def:mcarma} has a finite variance, it is also weakly stationary, i.e. the second-order moment structure (the variance and autocovariances) are left unchanged by time shifts.

Although the definition of CARMA processes via a spectral representation is elegant and helpful in many theoretical considerations, it is not really usable in applications, as alone simulating a CARMA process from this representation would be a tedious and problematic task. However, luckily we have the following result.
\begin{Theorem}[State Space Representation, \cite{Marquardtetal2005}]\label{th:statespace}
Let the L\'evy process $L$ and  $P,Q$ be  as before. Define the following
coefficient matrices:
\begin{itemize}
\item $\beta_{p-j}=-\sum\limits_{i=1}^{p-j-1} A_i\beta_{p-j-i}+B_{q-j}$, $j=0,1,\ldots,q$, $\beta_1=\ldots=\beta_{p-q-1}=0$
\item {${\beta}^*= \left(\beta_1^*,\beta_2^*,\ldots,\beta_p^*\right)$} and 
{${A}= \left(\begin{array}{c|ccc}0&&  I_{d(p-1)}& \\\hline
-A_p& -A_{p-1}&  \ldots& -A_1\end{array}\right)$}.
\end{itemize}
Denote by  {${G_t}=(G_{1,t}^*,\ldots,G_{p,t}^*)^*$ a
$pd$-dimensional process} and assume that  {$\mathcal{N}(P):=\{z\in\bbc:\det(P(z))=0\}\subset (-\infty,0)+i\bbr$} - the open right half of the complex plane. Then 
\begin{equation}\label{eq:sdeG}
dG_t=
{{A}
G_tdt+{\beta}dL_t}
\end{equation}
has a {unique stationary solution $G$} given by 
\begin{equation}\label{eq:G}
{
G_t=\int_{-\infty}^te^{{A}(t-s)}{\beta}\,dL_s},\quad t\in\bbr.
\end{equation}

It holds that \[ {G_{1,t}=\int_{-\infty}^{\infty}e^{i\mu t}P(i\mu)^{-1}Q(i\mu)\Phi_L(d\mu)=Y_t}, \:t\in\bbr.\]
So the {first $d$-components of $G$ are the CARMA process $Y$}.
\end{Theorem}

A CARMA process satisfying {$\mathcal{N}(P):=\{z\in\bbc:\det(P(z))=0\}\subset (-\infty,0)+i\bbr$} is called {causal}, 
because as shown above the value at a time $t$ only depends on the L\'evy process up to time $t$, it is a function of $(L_s)_{s\in(-\infty,t)}$. In other words a causal CARMA process is fully determined by values in the past. Whenever the condition $\mathcal{N}(P):=\{z\in\bbc:\det(P(z))=0\}\subset (-\infty,0)+i\bbr$ is not satisfied, $Y_t$ also depends on future values of the L\'evy process. In many applications, where it is clear that all we see today can only be influenced by what happened up to now, one only considers causal processes as appropriate models. However there are also applications where non-causal processes are useful. For example, if we want to stochastically model the water level in a river and think of $t$ as describing the location along the river, both the water levels downstream (in the ``future'') and upstream (in the ``past'') may influence the water level at a certain point.
Note that in this paper we only discuss stationary CARMA processes. In some applications (e.g. control) it is often adequate to consider non-stationary (non-stable) systems. Then the roots  $\det(P(z))=0\}$ in the set $(-\infty,0)+i\bbr$ describe the stable and causal part of the system and the remaining roots describe the non-stable part.

Theorem \ref{th:statespace} allows us to treat a causal CARMA process as a solution to the stochastic differential equation \eqref{eq:sdeG} and thus we can apply all the available results for SDEs. In particular, tasks like simulation of a causal CARMA process are straightforward and easily implemented. However, the above result allows us also to get rid of another restriction. So far we could only define CARMA processes driven by L\'evy processes with finite second moments and thus we could so far not have e.g. CARMA processes driven by $\alpha$-stable L\'evy processes. However, general theory on multidimensional Ornstein-Uhlenbeck processes (see \cite{Jureketal1993} and \cite{satoyamazato1984}) tells us that \eqref{eq:G} is the unique stationary solution to \eqref{eq:sdeG} as soon as the L\'evy process has only a finite logarithmic moment.
\begin{Definition}[Causal CARMA Process, \cite{Marquardtetal2005}]\label{def:causal}
Let $L=(L_t)_{t\in\bbr}$ be an $m$-dimen\-sional L\'evy process satisfying
\begin{equation}
\int\limits_{\|x\|\geq1}\ln\|x\|\,\nu(dx)<\infty,\label{lnu}
\end{equation}
$p,q\in\bbn_0$ with
$q<p$, and  further $A_1,A_2,\ldots,A_p,\in M_d(\bbr)$, $B_0,B_1,\ldots,\allowbreak B_q\in M_{d,m}(\bbr)$,
where  $B_0\not =0$. Define the matrices ${A}, {\beta}$ and the polynomial $P$ as in Theorem \ref{th:statespace} and assume  $\sigma(A)=\mathcal{N}(P)\subseteq (-\infty,0)+i\bbr$. Then the $d$-dimensional process 
\begin{equation}
Y_t=\left( I_d,0,\ldots,0\right)G_t
\end{equation}
where $G_t=\int_{-\infty}^te^{{A}(t-s)}{\beta} dL_s$ is the unique stationary solution to 
$
dG_t=
{A}G_tdt+
{\beta} dL_t
$
is called \emph{causal CARMA($p,q$)} process. 

$G$ is referred to as the state space representation.
\end{Definition}

A natural question is clearly whether one can also extend the definition of CARMA processes via the spectral representation to the case with infinite variance. For so-called  regularly varying L\'evy processes with finite mean and thus especially for $\alpha$-stable L\'evy processes with $\alpha\in (1,2)$ a result like Theorem \ref{th:rom} has been established in \cite{FuchsStelzer2010}. However, the non-finite variance case is distinctly different, as a limit of integrals has to be taken and the random orthogonal measure is replaced by an object which is -- strictly speaking -- not even a measure anymore. In that paper a definition of CARMA processes with regularly varying L\'evy input analogous to Definition \ref{def:mcarma} has been given and it has been shown that the resulting processes coincide with the causal CARMA processes when both definitions apply. Observe that processes with infinite variance are not only of academic interest, but that they have important applications, for instance, in network data modelling (cf. \cite{MikoschResnickRootzenStegeman2002} and \cite{Resnick2007,Resnick1997}). In \cite{GarciaKlueppelbergMueller2010} CARMA processes driven by $\alpha$-stable L\'evy processes have been successfully used to model electricity prices. 

\section{Properties}
In this section we explain and summarise various properties of (causal) CARMA processes.
\subsection{Second Order Structure}
Recall that for convenience we have assumed that the driving L\'evy process and thus the CARMA process has mean zero. Looking at the ``defining'' differential equations, it is clear that if $E(L_1)=\mu$ then the CARMA process is defined as the one driven by $L_1-\mu t$ plus $A_p^{-1}B_q\mu$  which is then the mean  of the CARMA process.

\begin{Proposition}[\cite{Marquardtetal2005}]
Let $Y$ be a (causal) CARMA process driven by a L\'evy process $L$ with finite second moments and set $\Sigma_L=\var(L_1)$.
\begin{enumerate}
\item The CARMA process $Y$ has autocovariance function: \[{\cov(Y_{t+h},Y_t)=\int\limits_{-\infty}^\infty
\frac{e^{i\mu h}}{2\pi}P(i\mu)^{-1}Q(i\mu)\Sigma_LQ(i\mu)^*(P(i\mu)^{-1})^*\,d\mu},\]
with $h\in\bbr.$
\item If $Y$ is a causal CARMA process, its state space representation $G$ has the following second order structure:
\begin{eqnarray*}
{\var(G_t)}&{=}&{\int\limits_0^\infty e^{{A}u}{\beta}\Sigma_L{\beta^*}e^{{A^*}u}du}\\
A{\var(G_t)}+{\var(G_t)}A^*&=&-{\beta}\Sigma_L{\beta^*}
\\{\cov(G_{t+h},G_t)}&{=}&{e^{{A}h}\var(G_t)}, \:h\geq 0.\end{eqnarray*}
\end{enumerate}
\end{Proposition}
Since we are only considering stationary CARMA processes, the moments above do not depend on $t$.

Since $Y$ is given by the first $d$ components of $G$ the second order structure of $G$ implies immediately alternative formulae for the second order structure of $Y$. In particular, it shows that the autocovariance function always decays like a matrix exponential for $h\to\infty$.
\subsection{Distribution}
Another nice feature is that in principle the distribution of a CARMA process at fixed times as well as the higher dimensional marginal distributions, e.g. the joint distribution of the process at two (or $n$) different points in time, is explicitly known in terms of the characteristic function. The reason is that all these distributions are infinitely divisible and that their L\'evy-Khintchine triplet  is known in terms of the L\'evy-Khintchine triplet of the driving L\'evy process. We state this in detail for the stationary distribution in the causal case.
\begin{Proposition} [\cite{Marquardtetal2005}]
If {$L$} has {characteristic triplet $(\gamma,\Sigma,\nu)$}, then the {stationary distribution} of the state space representation {$G$} of a {causal CARMA} process is {infinitely divisible} with {characteristic triplet $(\gamma_G^\infty,\Sigma_G^\infty,\nu_G^\infty)$}, where
\begin{eqnarray*}
&&\bullet\:\:{\gamma_G^\infty=
\int_0^\infty e^{{A}s}{\beta}\gamma\,ds+\int\limits_0^\infty\int\limits_{\bbr^{m}}
e^{{A}s}{\beta} x[1_{[0,1]}(\|{e^{{A}s}{\beta} x}\|)-1_{[0,1]}(\|{x}\|)]\,\nu(dx)\,ds,}\\
&&\bullet\:\:{\Sigma_G^\infty=\int_0^\infty e^{{A}s}{\beta}\,\Sigma{\beta}^*e^{{A}^*s}\,ds,} \\
&&\bullet\:\:{\nu_G^\infty(B)=\int_0^\infty\int_{\bbr^{m}}1_B(e^{{A}s}{\beta} x)\,\nu(dx)\,ds}
\end{eqnarray*}
for all Borel sets $B\subseteq\bbr^{pd}$.

In other words
\begin{align}
 E\left(e^{i\langle u,{G}_t\rangle}\right)&=\exp\left\{i\langle\gamma^\infty_G,u\rangle-\frac{1}{2}\langle u,\Sigma^\infty_G u\rangle+\int\limits_{\bbr^{pd}}
(e^{i\langle u,x\rangle}-1-i\langle u,x\rangle 1_{[0,1]}(\|x\|)\,\nu^\infty_G(dx)\right\},
\end{align}
for all $ u\in\bbr^{pd}$.
\end{Proposition}
Projection onto the first $d$ coordinates gives the characteristic triplet of the stationary distribution of  $Y$. It should, however, be noted that typically the distribution of the CARMA process does not belong to any special family of distributions even if one starts with especially nice L\'evy processes.

\subsection{Dependence Structure}
An important property of multivariate stochastic processes, is how their future evolution depends on the past. Suppose that one stands at a given point in time and one disposes of sufficient data at that point to determine the evolution from that point on, also given knowledge of the input from that point onwards. A {\it Markov process} is a stochastic process, for which the future only depends on the current value and not anymore on the past values (all their information is subsumed in the current value). For a Markov process it -- so to speak -- only matters where we are now not were we came from. If this characterising property does not only hold at all fixed times, but also at certain random times called stopping times, we speak of a strong Markov process.
\begin{Proposition} [\cite{Marquardtetal2005}]
The {state space representation $G$} of a {causal CARMA} process is a {strong Markov} process.
\end{Proposition}

Intuitively it is desirable in many applications that the farther away observations are in time, the less dependent they should be. Usually, one even wants that very far away observations should be basically independent. This idea is mathematically formalized in various concepts of asymptotic independence often referred to as some form of ``mixing''.

A comparably weak result which, however, applies to any CARMA process is the following.

\begin{Proposition}[\cite{FuchsStelzer2011}]
 Any stationary CARMA process is mixing. 
\end{Proposition}
Mixing implies ergodicity, i.e. empirically determined moments from the time series converge to the true moments if more and more data is collected. So time averages converge to  ensemble averages. This is very important for statistical estimation of CARMA processes, as it implies typically that estimators are consistent (i.e. the estimators converge to the correct value when more and more data is collected). 

Typically, one  also wants to know the errors of estimators which can be derived from distributional limit results like asymptotic normality. To obtain such results a stronger more uniform notion of asymptotic independence is needed, which is called strong mixing. Typically, one can best establish it for a Markov process.
\begin{Proposition}[\cite{Marquardtetal2005}]
 For a causal CARMA  process with $E(\|L_1\|^r)<\infty$ for some $r>0$ the state space representation  {$G$}  and the CARMA process {$Y$}  are {strongly mixing}, both with {exponentially decaying mixing coefficients}.
\end{Proposition}
\subsection{Sample Path Properties}
\begin{figure}
\begin{center}
	\includegraphics[width=0.8\textwidth]{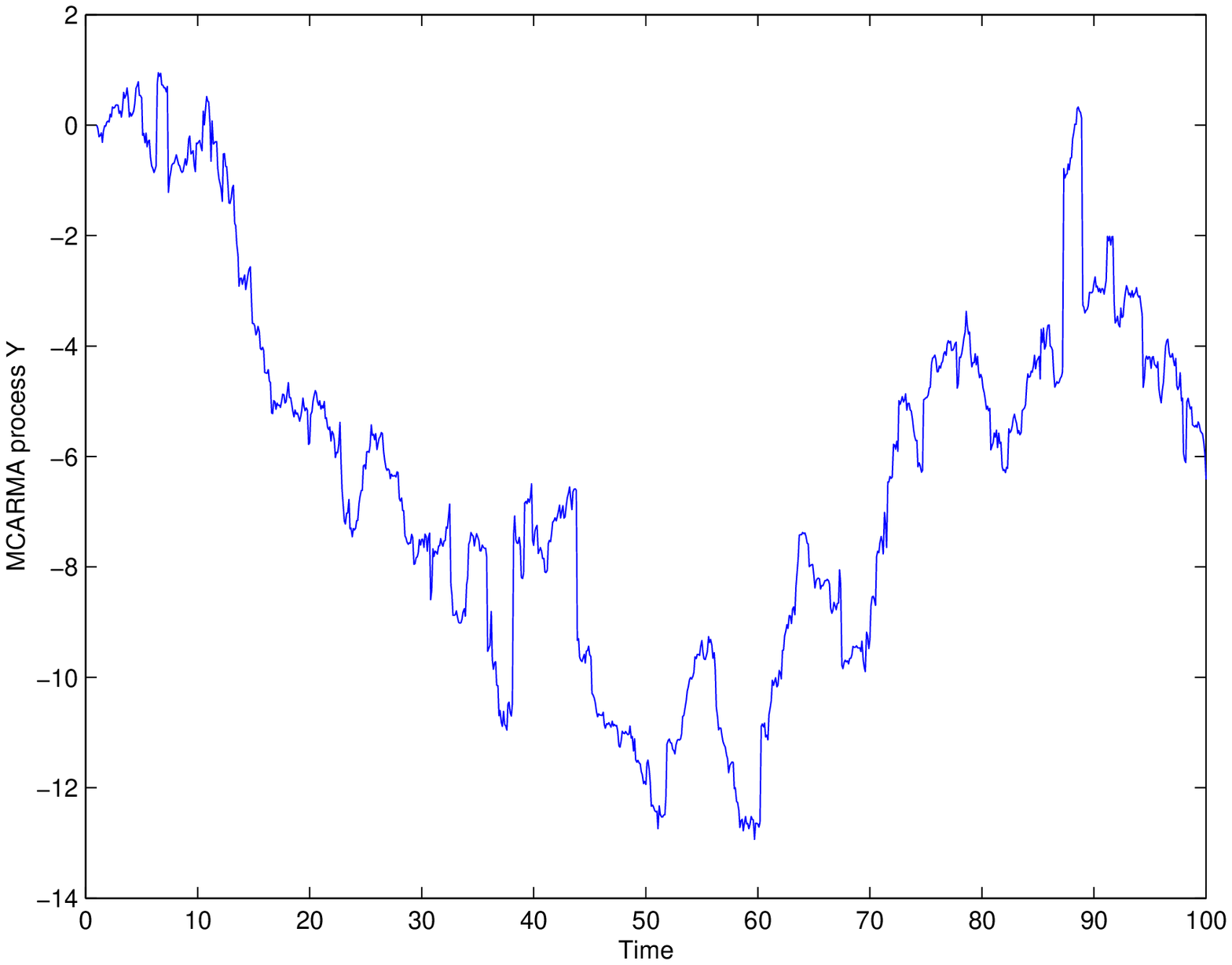}
	\includegraphics[width=0.8\textwidth]{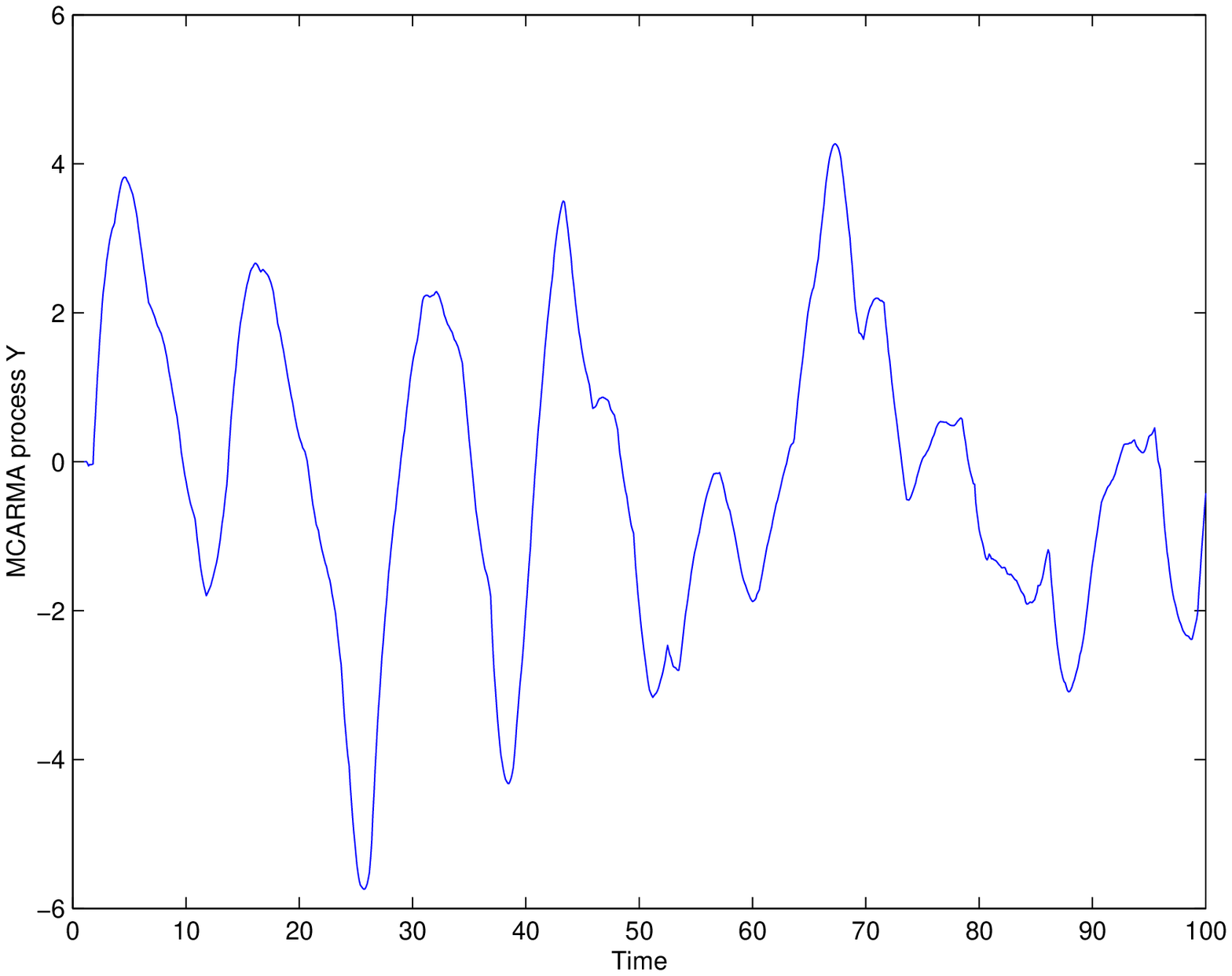}
\end{center}
\caption{A CARMA(1,0) process driven by an NIG L\'evy process having discontinuous paths is shown in the upper display and a CARMA(2,0) process driven by the same L\'evy process having continuous paths in the lower one. }\label{fig:carnig}
\end{figure}

Next we look at  the sample path properties of a CARMA process.
\begin{Proposition}[\cite{Marquardtetal2005}]\hfill
 \begin{itemize}
 \item The sample paths of a CARMA($p,q$) process {$Y$} with {$p>q+1$} are {($p-q-1$)-times 
differentiable} and for a causal CARMA process it holds that
\[\frac {d^i}{dt^i}Y_t=G_{i+1,t}, \:\:\:i=1,2,\ldots, p-q-1.\]
\item If {$p=q+1$} and the driving L\'evy process has a non-zero L\'evy measure $\nu$ satisfying $\nu(B_0^{-1}(\bbr^d\backslash\{0\}))\not=0$, then the paths of a CARMA process exhibit jumps and the jumps sizes are given by $\Delta Y_t:=Y_t-Y_{t-}=B_0\Delta L_t$.
\item If the driving L\'evy process {$L$   is a  Brownian motion}, 
then the sample paths of {$Y$} are {continuous} and {($p-q-1$)-times continuously differentiable}, provided $p>q+1$. 
\end{itemize}
\end{Proposition}
For examples of the paths of CARMA processes driven by an NIG L\'evy process see Figure \ref{fig:carnig}.

\subsection{Tail Behaviour}
As already stated in the introduction one may want to move away from Gaussian models, because extreme  (i.e. very low and/or high) observations are far more likely than in a Gaussian distribution. One says that the tails  (of the distribution) are heavier than Gaussian ones. Very often it appears also reasonable to use models which are ``heavy-tailed'' in the sense that only a limited number of moments exists, i.e. $E(\|X\|^r)$ exists only for low values of $r$. Mathematically it is then convenient to use the concept of regular variation (see \cite{Embrechtsetal1997} or \cite{Resnick2007,Resnick1987} for comprehensive introductions in relation to extreme value theory). Roughly speaking this means that the tails behave like a power function when one is far from the centre of the distribution. A random variable $X$ is regularly varying if $P(\|X\|>x)$ behaves comparably to $x^{-\alpha}$ for some $\alpha>0$ and big values of $x$. In \cite{Moseretal2010} (see also \cite{Fasen2005} in the univariate case) it is shown that under a very mild non-degeneracy condition a CARMA process driven by a regularly varying L\'evy process is again regularly varying with the same index $\alpha$. Hence, it is straightforward to construct heavy-tailed CARMA processes when applications call for such features.

In the univariate case the tail behaviour of CARMA processes is also understood in certain non-Gaussian situations, where one has lighter tails than regularly varying ones (see \cite{Fasen2009,fasen:2006}).

\section{State space models}

We have defined the causal CARMA process using a so-called state space representation and we have noted that the state space representation $G$ is made up of the CARMA process $Y$ and its derivatives as long as they exist. Hence, causal CARMA processes may be viewed as special state space models driven by L\'evy processes. In fact, any state space model can also be realized as a CARMA process, as will be shown now.

We start with a precise definition of state space models.
\begin{Definition}\label{Defstatespace}
Let $\Lb$ be an $m$-dimensional L\'evy process and \begin{equation*} A\in M_{N}(\R),\quad B\in M_{{N},m}(\R),\quad C\in M_{d,{N}}(\R).
 \end{equation*} A general $(N,d)$-dimensional continuous time state space model driven by $\Lb$ with parameters $A,B,C$ is a solution of 
\begin{align*}
\text{the state equation}\quad d\X_t =& A\X_tdt+Bd\Lb_t\\
\text{and the observation equation}\quad\Y_t =& C\X_t.
 \end{align*}

$\X$ is called the state process and $\Y$ the output process.
\end{Definition}

Note that the state  process is $N$-dimensional whereas the output process is $d$-dimensional.

Sufficient conditions for the existence of a unique causal stationary solution of the state equation are given by ($\Re(\cdot)$ indicates the ``real part'' of a complex number or function)
 \begin{equation*}
 \Re(\lambda_\nu)<0,\quad\lambda_\nu,\,\nu=1,\ldots,{N}, \text{ being the eigenvalues of } A
 \end{equation*}
 and $\Lb$ having  finite second moments.

It can easily be shown by integration that $\X$ satisfies
\begin{equation*}
\X_t = e^{A(t-s)}\X_s+\int_s^t{e^{A(t-u)}Bd\Lb_u}.
\end{equation*}
Likewise, the stationary {output process $\Y$} satisfies 
\begin{equation*}
\Y_t = \int_{-\infty}^t{Ce^{A(t-u)}Bd\Lb_u}.
\end{equation*}
Its {spectral density}, the Fourier transform of the autocovariance function, is given by 
{\begin{equation*}
f_{\Y}(\omega) = \frac{1}{2\pi}C(i\omega-A)^{-1}B\Sigma_{\Lb}B^T(-i\omega-A^T)^{-1}C^T.
\end{equation*}}

From Definition \ref{def:causal} it is obvious that a CARMA process is a $(pd,d)$-dimensional state space model driven by an $m$-dimensional L\'evy process. The following theorem states that also the converse is true.
\begin{Theorem}[\cite{Schlemm:Stelzer:2010a}]\label{th:statespacecarma}
The stationary solution $\Y$ of the multivariate state space model $(A,B,C,\Lb)$ is an $\Lb$-driven CARMA process with autoregressive polynomial $P$ and moving average polynomial $Q$ if and only if
\begin{equation*}
C(zI_{N}-A)^{-1}B = P(z)^{-1}Q(z),\quad\forall z\in\mathbb{C}.
\end{equation*}

For any $(A,B,C)$ there exist $P,Q$ such that the above equation is satisfied and vice versa.
\end{Theorem}

In reality we typically do not observe some variables of interest continuously, but only at a discrete set of points in time. Let us assume that we sample the process at an equidistant time grid with grid  length $h>0$ and denote by $\Y_n^{(h)}:=\Y_{nh}$ for   $n\in\bbz$ the sampled observations of a state space process.

It is easy to see that
\begin{align}
 \Y_n^{(h)}&=C\X_n^{(h)}\\
 \X_n^{(h)}&= e^{Ah}\X_n^{(h)}+\int_{(n-1)h}^{nh}e^{A(nh-u)}Bd\Lb_u,
\end{align}
which immediately shows that $\Y_n^{(h)}$ is the output process of a discrete time $(N,d)$-dimensional state space model driven by the $N$-dimensional iid noise $\left( \int_{(n-1)h}^{nh}{e^{A(nh-u)}}Bd\Lb_u\right)_{n\in\bbz}$.

It is well-known that any $(N,d)$-dimensional state space model in discrete time is an ARMA process. Combining this with Theorem \ref{th:statespacecarma} tells us that any equidistantly sampled CARMA process $Y^{(h)}$ is an ARMA process. This observation will be the basis for estimating CARMA parameters in the next section, where we will need a considerable refinement of this result.

In many applications the sampling frequency is quite high, i.e. $h$ is very small. Thus it is important to understand how $Y^{(h)}$ behaves as $h\to 0$ which has been investigated in \cite{BrockwellFerrazzanoKlueppelberg2011}.

As we only observe the process $\Y$ in a state space model, an important question is what can be said about the state process $\X$ based on the observations. Hence, we want to reconstruct or ``estimate'' the latent process $\X$ as good as possible. This procedure is also referred to as filtering. For Gaussian state space models the easily implementable Kalman filter (see e.g. \cite{Brockwelletal1991}) is optimal both from a variance point as well as a distributional point of view. For non-Gaussian state space models with finite variance the very same procedure, now typically called linear filtering, gives an ``estimate'' of the latent process which is the linear (in the observations) ``estimate'' with the lowest variance. However, it is typically not the ``estimate'' with the minimal variance and not a conditional expectation. Thus, there are more involved filtering techniques, like particle filtering (see e.g. \cite{DoucetFreitasGordon2001}), which are better.

State space models, mainly Gaussian ones, are also heavily used in stochastic control (see \cite{GarnierWang2008} and references therein for a comprehensive overview) and signal processing (see \cite{Mossbergetal2004,LarssonMossbergSoederstroem2006,Larsson2005}, for instance). In both areas one is sometimes dealing with data for which a Gaussianity assumption is not really appropriate due to skewedness, excess kurtosis or heavy-tailedness. Clearly, in such situations L\'evy-driven state space models or equivalently CARMA processes should be appealing. 
Going into the details of the usage in control is beyond the scope of this paper, but it seems worthwhile to mention that there are two uses of state space models in control. Sometimes one  assumes that one has some random input which is then ``controlled'' by the state space model, so the the state space model acts as the controller. In contrast to this sometimes the output of the state space model is regarded as the natural output of some system on which an additional controller is acting to ensure that the output meets certain requirements.

\section{Statistical Estimation}
In this section we discuss ways to estimate the parameters of a CARMA process and its driving L\'evy process. First we address the estimation of the autoregressive and moving average parameters. Due to parametrisation issues explained later on, we formally do this for L\'evy-driven continuous time state space models, as defined in the previous section. In the univariate case quasi-maximum likelihood estimation of CARMA processes is comprehensively studied in \cite{BrockwellDavisYang2009}.
\subsection{Quasi-maximum likelihood estimation}
We assume that we observe the process $\Y$ at discrete, equally spaced times
\begin{equation*}
\Y_n^{(h)}:= \Y_{nh},\quad n\in\mathbb{Z},\quad h>0.
\end{equation*}
Furthermore, we define the linear innovations ${\varepsilon}^{(h)}$ by
\begin{align*}
\boldsymbol{\varepsilon}^{(h)}_n = \Y^{(h)}_n-&P_{n-1}\Y^{(h)}_n,\quad 
\end{align*}
where $P_{n-1}$ denotes the orthogonal projection onto $\overline{\operatorname{span}}\left\{\Y^{(h)}_\nu:-\infty<\nu<n\right\}$, i.e. the linear space spanned by the observations until time $(n-1)h$.
From the construction it is immediate that $(\boldsymbol{\varepsilon}^{(h)}_n)_{n\in\mathbb{Z}}$ is a white noise sequence, i.e. it has mean zero, a constant variance and is uncorrelated. The construction implies that one can only sensibly speak of linear innovations when the driving L\'evy process has finite second moments. Thus we will  demand the latter for the remainder of this section.

\begin{Theorem}[\cite{Schlemm:Stelzer:2010a}] Assume the eigenvalues $\lambda_1,\ldots, \lambda_N$ of the matrix $A$ are pairwise distinct and
define complex numbers $\Phi_1, \Phi_2,\ldots, \Phi_N$ by
\begin{equation*}
1-\Phi_1z-\Phi_2 z^2-\ldots-\Phi_Nz^N = \prod_{\nu=1}^{{N}}\left[1-e^{-\lambda_\nu h}z\right]\,\forall\, z\in \bbc.
\end{equation*}
Then there exist $\Theta_1,\Theta_2,\ldots,\Theta_{N-1}$ in $M_d(\mathbb{C})$  such that
\begin{equation*}
\Y_n^{(h)}-\Phi_1\Y_{n-1}^{(h)}-\ldots-\Phi_N\Y_{n-N}^{(h)}={\varepsilon}^{(h)}_n+\Theta_1{\varepsilon}^{(h)}_{n-1}+\ldots+\Theta_{N-1}{\varepsilon}^{(h)}_{n-N+1}
\end{equation*}
holds.

Hence, $\Y^{(h)}$ is a weak $\operatorname{ARMA}({N},{N}-1)$ process.
\end{Theorem}
This result suggests that one could estimate simply the ARMA coefficients of the sampled process and then transfer these estimates to estimates of the CARMA coefficients.
However,
to estimate a CARMA process it is not sufficient to estimate an ARMA process, because not all ARMA processes can be embedded in a CARMA process. There are ARMA processes which cannot arise as equidistantly sampled CARMA processes. The way out is carry out the ``ARMA estimation''  in  the CARMA parameter space.

Since we are going to use a quasi-maximum likelihood approach  and have discretely sampled observations, all possible models considered in the estimation have to be distinguishable  based only on the {second-order properties of the sampled process}.

\begin{Definition}[Identifiability]
A collection of continuous time stochastic processes $\left(\Y_{\bth},\bth\in\Theta\right)$ is {identifiable} if for any $\bth_1\neq\bth_2$ the two processes $\Y_{\bth_1}$ and $\Y_{\bth_2}$ have different spectral densities.\

It is {$h$-identifiable}, $h>0$, if for any $\bth_1\neq\bth_2$ the two processes $\Y_{\bth_1}^{(h)}$ and $\Y_{\bth_2}^{(h)}$ have different spectral densities.
\end{Definition}

We assume that our parametrisation is given by a compact parameter space $\Theta\subset \bbr^q$ with some $q\in\bbn$  and a mapping 
\begin{equation*}
\psi:\Theta\ni\bth\mapsto(A_{\bth},B_{\bth},C_{\bth},\Lb_{\bth}).
\end{equation*}
Here, $A_\vartheta$ is the $N\times N$ matrix of our Definition \ref{Defstatespace} dependent on the parameters $\vartheta$ and likewise for $B_{\bth},C_{\bth}$ and $\Lb_{\bth}$.
 
We need to ensure that our parametrisation is minimal regarding the dimensions, since a fixed output process can result from artificially  arbitrarily high-dimensional state space models.
\begin{assumptionpara}[Minimality]
For all $\bth\in\Theta$ the triple $\left(A_{\bth} ,B_{\bth} ,C_{\bth} \right)$ is {minimal} in the sense that if
\begin{equation*}
C(zI_{m}-A)^{-1}B = C_{\bth}(zI_{N}-A_{\bth})^{-1}B_{\bth}
\end{equation*} then $m\geq N$ must be true.
\end{assumptionpara}
\begin{assumptionpara}[Eigenvalues]
For all $\bth\in\Theta$ the {eigenvalues of $A_{\bth}$} are pairwise distinct and contained in the strip 

\begin{equation*}
\{z\in\mathbb{C}:-\pi/h<\Im( z) < \pi/h\}.
\end{equation*}
\end{assumptionpara}

We want to use a parametrisation for the continuous time state space model, but need to ensure that it is $h$-identifiable. The following theorem provides easy-to-check criteria.
\begin{Theorem}[\cite{Schlemm:Stelzer:2010b}]
Assume that the parametrisation $\psi:\Theta\supset\bth\mapsto\left(A_{\bth} ,B_{\bth} ,C_{\bth} ,\Lb_{\bth}\right)$ is 
\begin{itemize}
 \item identifiable,
 \item minimal
 \item and satisfies the eigenvalue condition.
\end{itemize}
 Then the corresponding collection of output processes $\left\{\Y_{\bth},\bth\in\Theta\right\}$ is {$h$-identifiable}.
\end{Theorem}

The quasi-maximum likelihood (QML) estimator is now obtained by pretending the observations were Gaussian, taking the corresponding likelihood and maximising it. More precisely the QML of  $\bth$ based on $L$ observations $\y^L=(\y_1,\ldots,\y_L)$ (of a CARMA process with parameter $\vartheta_0$) is
\begin{equation*}
{\hat{\bth}^L= \operatorname{argmax}_{{\bth}\in\Theta}{\mathcal{L}_{\bth}\left(\y^L\right)}
},
\end{equation*}
where $\mathcal{L}_{\bth}$ is the {Gaussian likelihood function} which is proportional to
 \begin{equation*}
\left(\prod_{n=1}^L{\det V_{\bth,n}}\right)^{-1/2}\exp\left\{-\frac{1}{2}\sum_{n=1}^L{\mathbf{e}_{\bth,n}^TV_{\bth,n}^{-1}\mathbf{e}_{\bth,n}}\right\}
 \end{equation*}
with
 \begin{align*}
 \mathbf{e}_{\bth,n}=&\y_n-\left.P_{n-1}\Y_{\bth,n}^{(h)}\right|_{\left\{\Y_{\bth,\nu}^{(h)}=\y_\nu:1\leq\nu<n\right\}},\\
 V_{\bth,n}=&\mathbb{E}\left[\mathbf{e}_{\bth,n}\mathbf{e}_{\bth,n}^T\left|\Y_{\bth,\nu}^{(h)}=\y_\nu:1\leq\nu<n\right.\right].
\end{align*}
So $\mathbf{e}_{\bth,n}$ are the linear innovations under the model given by $\vartheta$ and $V_{\bth,n}$ are their variances or the one-step prediction errors. Note that $\y^L$ are in contrast to this observations of the CARMA process with the unknown parameter $\vartheta_0$ which we are about to estimate.

Computing the QML estimator is now a straightforward task utilising the Kalman recursions and numerically maximising the likelihood. However, since we have not used the true likelihood, it is not clear whether the resulting estimators are really sensible in the sense that they converge to the true parameters. Luckily, one can show that the estimators are well-behaved.

\begin{Theorem}[Strong consistency, \cite{Schlemm:Stelzer:2010b}]
Assume the parametrisation $\psi$ is continuous. For every sampling interval $h>0$, the QML estimator $\hat\bth^L$ is strongly consistent, i.e.
\begin{equation*}
\hat{\bth}^L\to\bth_0\quad \text{a.s. as } L\to\infty,
\end{equation*}
provided the parametrisation is $h$-identifiable.
\end{Theorem}

However, so far we cannot assess the quality of our estimators by confidence intervals etc., which is made possible by the following result.

\begin{Theorem}[Asymptotic normality, \cite{Schlemm:Stelzer:2010b}]\label{th:asymnorm}
Assume that the driving L\'evy process satisfies ${E}||\Lb_{\vartheta_0}(1)||^{4+\delta}<\infty$ for some $\delta>0$ and that the parametrisation $\psi$ is three times continuously differentiable. For every sampling interval $h>0$, the QML estimator $\hat\bth^L$ is asymptotically normally distributed, i.e.
\begin{equation*}
\sqrt{L}\left(\hat{\bth}^L-\bth_0\right)\stackrel{\mathcal{D}}{\to}\mathcal{N}(0,\Omega),\quad \Omega=J(\bth_0)^{-1}I(\bth_0)J(\bth_0)^{-1}, 
\end{equation*}
with
\begin{align*}
J(\bth) =&\lim_{L\to\infty}{\frac{1}{L}\frac{\partial^2}{\partial\bth\partial\bth^T}\ln \mathcal{L}_{\bth}\left(\y^L\right)},\\
I(\bth)=&\lim_{L\to\infty}{\frac{1}{L}\operatorname{Var}{\frac{\partial}{\partial\bth}\ln \mathcal{L}_{\bth}\left(\y^L	\right)}},
\end{align*} provided the parametrisation is $h$-identifiable.
\end{Theorem}

To obtain identifiable parametrisations one uses like in the discrete time case (see \cite{hannan1987stl} or \cite{Luetkepohl2005}, for instance) so called canonical parametrisations like the echelon state space form. For more details on this we refer to \cite{Schlemm:Stelzer:2010b}. Since such parametrisations are typically available for state space models rather than CARMA processes, one normally estimates state space models rather than the equivalent CARMA processes.

Let us finally look at one simulation study.

A $d$-dimensional normal inverse Gaussian (NIG)  L\'evy process $\Lb$ (see e.g. \cite{Prause1999,Blaesildetal1981,Barndorff1998}) with parameters
\begin{equation*}
\delta>0,\kappa>0,\boldsymbol{\beta}\in\R^d,\Delta\in M_d^+(\R)
\end{equation*}
is given by a normal mean-variance mixture, i.e.
\begin{equation*}
\Lb_1{=}\boldsymbol{\mu} + V\Delta\boldsymbol{\beta}+V^{1/2}N,
\end{equation*}
where $N$ is $d$-dimensionally  normally distributed with mean zero and variance $\Delta$ and independent of  
\begin{equation*}
V\sim\operatorname{IG}(\delta/\kappa,\delta^2)
\end{equation*}
which follows a so-called inverse Gaussian distribution (\cite{Joergensen1982}). 

We consider now a
{bivariate} NIG-driven CARMA process with zero mean given by the state space form

\begin{align*}
d\X_t=&\left[\begin{array}{ccc}
		\vartheta_1 & \vartheta_2 & 0 \\
		0 & 0 & 1 \\
		\vartheta_3 & \vartheta_4 & \vartheta_5
	    \end{array}\right]\X_tdt+\left[\begin{array}{cc}
					    \vartheta_1 & \vartheta_2 \\
					    \vartheta_6 & \vartheta_7 \\
					    \vartheta_3+\vartheta_5\vartheta_6 & \vartheta_4+\vartheta_5\vartheta_7
					\end{array}\right]d\Lb_t,\\
\Y_t=&\left[\begin{array}{ccc}1 & 0 & 0 \\ 0 & 1 & 0\end{array}\right]\X_t,\quad \Sigma_{\Lb} =\left[\begin{array}{cc} \vartheta_8 & \vartheta_9 \\\vartheta_9 & \vartheta_{10}\end{array}\right].
\end{align*}
The parameters are $\vartheta_1,\vartheta_2,\ldots,\vartheta_{10}$ and the parametrisation is in one of the canonical identifiable forms.

\begin{figure}[tp]
\begin{center}
\includegraphics[width=0.8\textwidth]{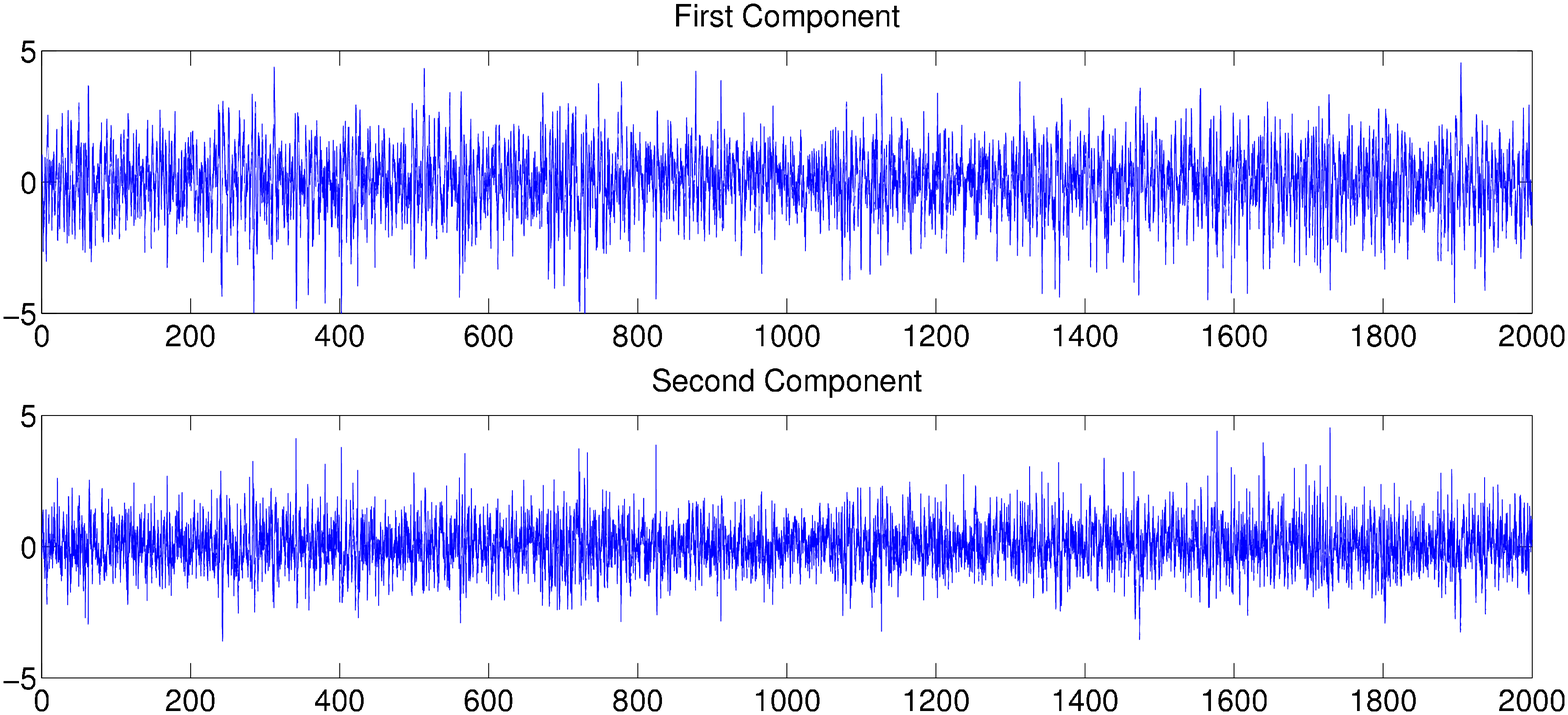}\vspace*{0.5cm}

\includegraphics[width=0.8\textwidth]{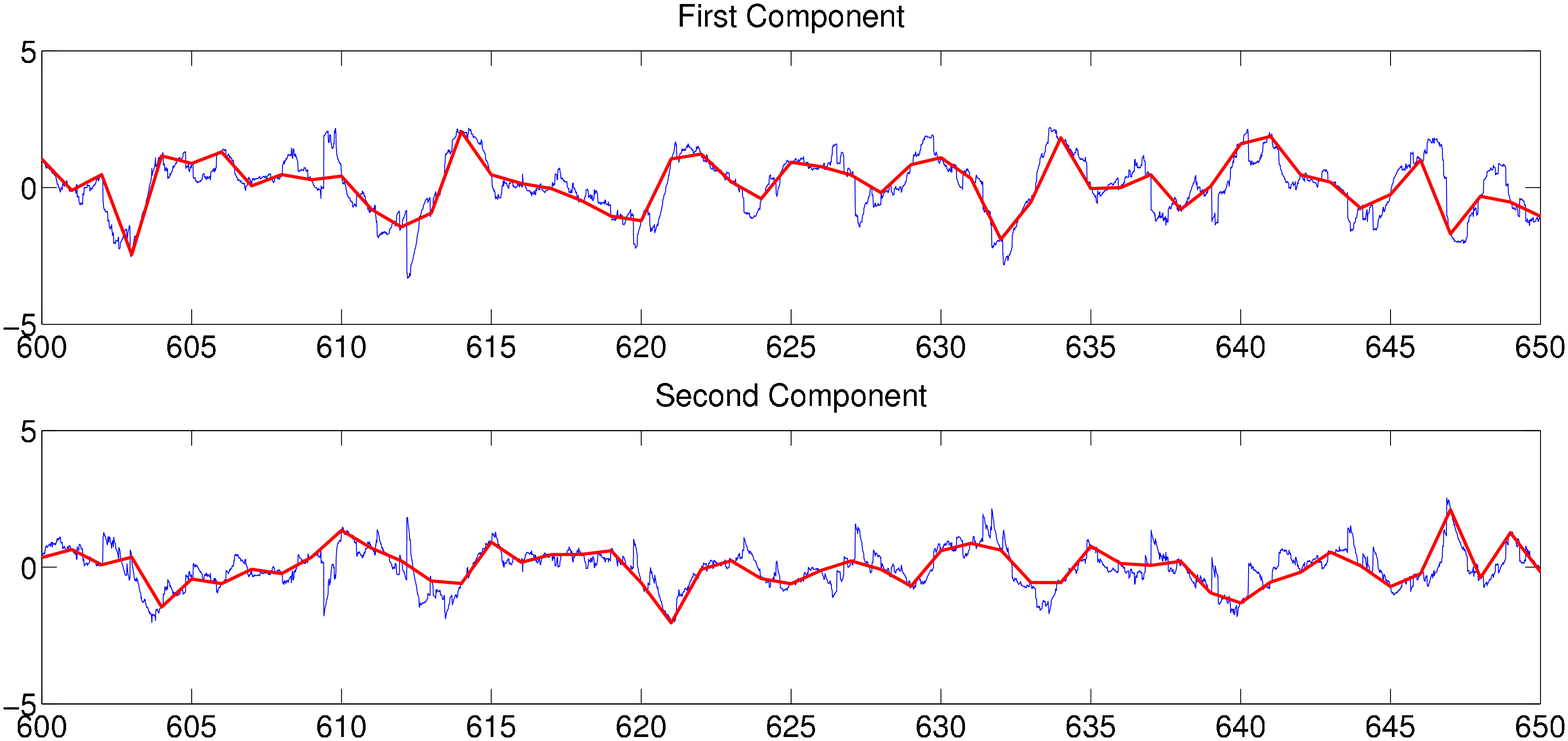}
\end{center}
\caption{One {realisation} of a bivariate NIG-driven CARMA process (upper two displays) and the effect of sampling (lower two displays). The linearly interpolated  process over the time interval $[600,650]$ resulting from sampling at integer times is shown as the thicker line, whereas the thinner line is the true CARMA process.}\label{fig:sim}
\end{figure}

A simulated path is shown in Figure \ref{fig:sim}.

We calculated the QML estimates for this {bivariate NIG-driven CARMA process based on 
observations over the time horizon $[0,2000]$ at integer times and repeated this for 350 different simulated paths. The estimation results are summarised in Table \ref{table:est}. It shows that the sample bias of the obtained estimators in the simulation study is very small and that the sample standard deviation is close to the standard deviation predicted by the asymptotic normality result Theorem \ref{th:asymnorm}. Actually, the sample standard deviation is always smaller which is nice, as it implies that the standard deviation predicted by the asymptotic normality result Theorem \ref{th:asymnorm} is a conservative estimate.
\begin{table}
\centering
\begin{tabular}{|c||c|c|c|c|}
 \hline
parameter & sample mean & sample bias & sample & estimated \\
 & & & standard deviation & standard deviation\\
\hline\hline
$\vartheta_1$ 	& -1.0001	& 0.0001	&0.0354		&0.0381	\\
$\vartheta_2$ 	&  -2.0078	& 0.0078	&0.0479		&0.0539	\\
$\vartheta_3$ 	&  1.0051	&-0.0051	&0.1276		&0.1321	\\
$\vartheta_4$ 	&  -2.0068	& 0.0068	&0.1009		&0.1202	\\
$\vartheta_5$ 	&  -2.9988	&-0.0012	&0.1587		&0.1820	\\
$\vartheta_6$ 	&  1.0255	&-0.0255	&0.1285		&0.1382	\\
$\vartheta_7$ 	&  2.0023	&-0.0023	&0.0987		&0.1061	\\
$\vartheta_8$ 	&  0.4723	&-0.0028	&0.0457		&0.0517	\\
$\vartheta_9$ 	&  -0.1654	& 0.0032	&0.0306		&0.0346	\\
$\vartheta_{10}$&  0.3732	& 0.0024	&0.0286		&0.0378	\\
\hline
\end{tabular}
\caption{Summary of the results of the simulation study on the QML estimation of a bivariate NIG-driven CARMA process. The second column states the mean of estimators obtained over 350 simulated paths, the third column the resulting bias and the fourth column the standard deviation of the obtained estimators. Finally, the last column states the standard deviation for the estimators as predicted by the asymptotic normality result Theorem \ref{th:asymnorm}.}\label{table:est}
\end{table}

\subsection{Statistical inference for the driving L\'evy process}\label{sec:levyest}
The above quasi-maximum likelihood approach only allows to estimate the autoregressive and moving average parameters as well as the variance of the driving L\'evy process. However, typically we want to estimate many more parameters of the driving L\'evy process or even first need to get an idea to which family the driving L\'evy process may belong to. To this end one can reconstruct from the CARMA process the driving L\'evy process. Typically, the CARMA process is only observed at a discrete set of times and then the best we can do is to get approximations of the increments of the L\'evy process. One can then treat the approximate increments as if they were the true ones of the L\'evy process. ``Looking'' at them one should be able to choose appropriate parametric families. By using the approximate increments, as one would use the true ones, in maximum likelihood or method of moment based estimation procedures one can  do parametric inference for the L\'evy process. The construction of the approximate increments and their use in estimation procedures has been studied in detail in \cite{BrockwellSchlemm2011} where it is in particular shown that the estimators are good  in the sense that they are consistent and asymptotically normal under reasonable assumptions when taking appropriate limits.

It should be noted that the idea to reconstruct the L\'evy process can already be found in \cite{Brockwell2009} or \cite{BrockwellDavisYang2009}. In the following we illustrate this approach for a univariate Ornstein-Uhlenbeck, i.e. a CARMA(1,0), process based on \cite{BrockwellDavisYang2007b} and \cite{Graf2009} from which all examples and plots are taken.

Recall that an Ornstein-Uhlenbeck (OU) process is the unique strictly stationary solution to
\begin{equation} \label{DefOU}
dY_t=a Y_tdt+ dL_t.
\end{equation}  
where  $\left(L_t\right)_{t\in\bbr}$ is a L\'evy process with $E(\ln(\max(|L_1|,1)))<\infty$ and autoregressive parameter $a < 0$.
The solution of the stochastic differential equation is given explicitly by
\begin{equation}\label{ProOU} Y_t=e^{a(t-s)}Y_s+ \int\limits_{s}^{t}{ e^{a(t-u)}dL_u}.
\end{equation}

If the OU process is observed continuously on $\left[0,T\right]$, then the integrated form of \eqref{DefOU} immediately gives
\begin{equation*}
L_t=Y_t-Y_0-a\int_{0}^{t}{Y_s ds}.
\end{equation*}
The increments of the driving L\'evy process $\Delta L_n^{(h)}$ on the intervals $\left((n-1)h,nh\right]$ with $n\in\bbn$ can be represented as
\begin{equation}\label{increL}
\Delta L_n^{(h)}:=L_{nh}-L_{(n-1)h}=Y_{nh}-Y_{(n-1)h}-a\int_{(n-1)h}^{nh}{Y_u du}.
\end{equation}

What we want, is to approximately reconstruct the sequence $\Delta L^{(h)}$ of increments over intervals of length $h$ from observations of the CARMA process made over a finer equidistant grid. To this end one simply approximates the integral $\int_{(n-1)h}^{nh}{Y_u du}$ by some numerical integration scheme needing only the values of the process on this finer grid.
Since the approximations of $\Delta L^{(h)}$ become thus closer and closer to the true increment as the numerical integration scheme becomes more exact,  \cite{BrockwellSchlemm2011} derive their asymptotic results when both the observation interval as well as the observation frequency goes to infinity. 
Note that in practice one does not know $a$ so one has to estimate it first, which could e.g. be done by the already described quasi-maximum likelihood approach.

Turning to an example, let us consider the OU process given by 
\begin{equation}
dX_t=-0.6 X_tdt+ dL_t,
\end{equation}
with $L$ being a standardised Gamma process, i.e. $L_t$ has density
\begin{equation*}
f_{L_t}(x)=\frac{\gamma^{1/2 \gamma t}}{\Gamma(\gamma t)}x^{\gamma t-1}e^{-x \gamma^{1/2}}\textbf{1}_{\left[0,\infty\right)},
\end{equation*}
and the parameter $\gamma$ being set to 2.

In \cite{Graf2009}  100 paths of this OU process on the time interval $[0,5000]$ have been simulated and then the L\'evy increments over time intervals of unit length have been approximated by sampling the OU process over a grid of size $h$.

In Figure \ref{fig:EulerCompTrap} the histogram of the L\'evy increments distribution from one path with $h=0.01$ is shown, together with the true probability density of $L_1$. 
\begin{figure}[tp]
\centering
\includegraphics[width=0.7\textwidth]{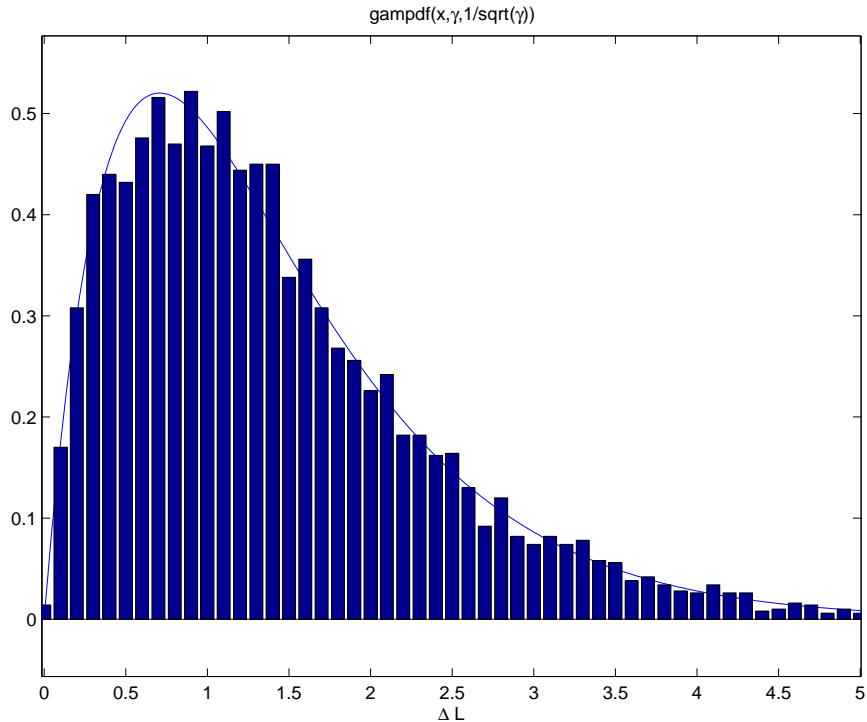}
\caption{Probability density of the increments of the standardised L\'evy process with $\gamma =2$ and the histogram of the estimated increments from one path of the OU process, obtained by sampling the process with grid length 0.01. (Source: \cite{Graf2009}.)}\label{fig:EulerCompTrap}
\end{figure}

If one further averages over all one hundred paths which is equivalent to looking at one path over a one hundred times longer time horizon, the fit of the histogram to the true density becomes visually almost perfect, see Figure \ref{fig:EulerCompTrapAll}.
\begin{figure}[tp]
\centering
\includegraphics[width=0.7\textwidth]{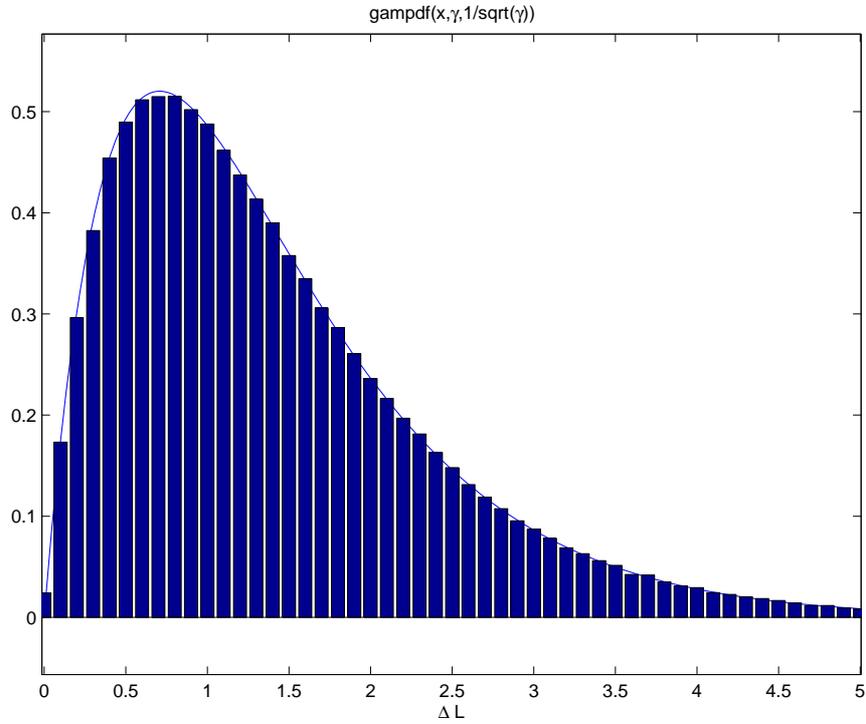}
\caption{Probability density of the increments of the standardised L\'evy process with $\gamma =2$ and the histogram of the estimated increments for all 100 paths of the OU process, obtained by sampling the process with grid length 0.01. (Source: \cite{Graf2009}.)}\label{fig:EulerCompTrapAll}
\end{figure}

Based on the approximate L\'evy increments one can now estimate the parameter $\gamma$ by maximum likelihood.
Table \ref{TabDATAGamma}  shows summary statistics of the resulting  estimator for different samp\-ling grid sizes $h$. The data in the table is based on estimating $\gamma$ separately for each of the 100 simulated paths. 

\begin{table}[tp]
\caption{Estimated parameters of the standardised driving L\'evy process based on 100 paths on $\left[0,5000\right]$ of the Gamma-driven OU process.}
\begin{center}
\begin{tabular}{|c||c|c|c|} \hline
$h$ & Parameter & Sample mean & Sample standard\\
 & &     of estimator & deviation of estimator\\ \hline\hline
0.01 & $\gamma$ & 2.0039 & 0.0314\\
 0.1& $\gamma$& 2.0043 & 0.0340\\
1 & $\gamma$ & 1.9967 & 0.0539\\ \hline
\end{tabular}
\end{center}
\label{TabDATAGamma}
\end{table}

To conclude, the simulation study illustrates that the recovery of the background driving L\'evy process and the parametric estimation based on the approximate increments works quite well.

\section{Concluding Remarks}
Finally, we would like to mention that there are other stochastic models like the so-called ECOGARCH process of \cite{Haugetal2006} and \cite{HaugetStelzer2008} where CARMA processes are an important ingredient as well as extensions of CARMA processes. One extension are fractionally integrated CARMA  (FICARMA) processes (see \cite{Brockwelletal2005} and \cite{Marquardt2007}).  While CARMA processes have an exponentially decaying autocovariance function and thus have always short memory, FICARMA processes exhibit polynomially decaying autocovariance functions and are thus able to model long memory phenomena (see \cite{Doukhanetal2003} or \cite{Samorodnitsky2006} for detailed introductions into the topic of long range dependence). However, the paths of FICARMA processes are continuous. A class of processes with possible long memory, jumps in the paths and related to CARMA processes are the supOU processes, see \cite{BarndorffStelzer2009,barndorffnielsen01} and \cite{FasenetCklu2007}. As noted in \cite{BarndorffStelzer2009} multivariate supOU processes can be straightforwardly extended to obtain so-called supCARMA processes.
Long memory is (believed to be) encountered in data from many different areas, e.g. finance or telecommunication. Since it is an asymptotic property and similar effects in the autocorrelation function might be caused by structural breaks (non-stationarity), it is often hardly debated whether there truly is long memory in a time series. The first scientific study considering long range dependence properties was looking at the water level of the river Nile (see \cite{Tousson.1925}).

From the overview on  CARMA processes presented in this paper it should not only be clear that they are useful in many applications, but also that there are still many questions to be addressed in future research. These include alternative estimators to the ones presented here, estimators which work in the heavy-tailed case when one does not have a finite variance or order selection, i.e. a theory how to choose the orders $(p,q)$ of the autoregressive and moving average polynomial when one fits CARMA processes to observed time series.

\section*{Acknowledgements}
The author takes pleasure in thanking Florian Fuchs, Claudia Klüppelberg and Eckhard Schlemm for comments on previous drafts and Maria Graf for allowing him to use material from her diploma thesis in Section \ref{sec:levyest}.
Financial support  from the TUM Institute for Advanced Study funded by the German Excellence Initiative through a Carl-von-Linde Junior Fellowship is gratefully acknowledged.

\end{document}